\documentclass{scrartcl}

\usepackage[utf8]{inputenc}
\usepackage[USenglish]{babel}

\usepackage{pgfplots}
\usetikzlibrary{arrows,decorations,decorations.pathreplacing}
\usetikzlibrary{plotmarks}
\pgfplotsset{compat=newest}
\usepackage{mathtools,amsmath,amsthm,amssymb,esint}
\usepackage[unicode=true,pdfusetitle,
 bookmarks=true,bookmarksnumbered=false,bookmarksopen=false,
 breaklinks=true,pdfborder={0 0 1},backref=false,colorlinks=true]
 {hyperref}
\usepackage{subfig,cleveref}

\theoremstyle{plain}
  
  \newtheorem{prop}{Proposition}
  
\theoremstyle{remark}
  \newtheorem{rem}{Remark}
\theoremstyle{definition}
  \newtheorem*{definition}{Definition}
  \newtheorem{example}{Example}
  
\DeclareMathOperator{\linhull}{span}
\DeclareMathOperator{\grad}{grad}
\DeclareMathOperator{\divergence}{div}
\DeclareMathOperator{\id}{Id}

\DeclareMathOperator{\dist}{dist}
\DeclareMathOperator{\Vol}{Vol}

\DeclareMathOperator{\sech}{sech}
\DeclareMathOperator{\trace}{trace}

\title{A Geometric Heat-Flow Theory of\\ Lagrangian Coherent Structures}
\author{Daniel Karrasch\thanks{Electronic mail: \href{mailto:karrasch@ma.tum.de}{karrasch@ma.tum.de}}%
\and Johannes Keller \and \\%
Technische Universität München, Zentrum Mathematik\\
Boltzmannstr.~3, 85748 Garching bei München, Germany}
\date{\today}

\begin{document}

\maketitle

\begin{abstract}
We consider Lagrangian coherent structures (LCSs) as the boundaries of material
subsets whose advective evolution is metastable under weak diffusion.
For their detection, we first transform the Eulerian advection--diffusion equation
to Lagrangian coordinates, in which it takes the form of a time-dependent diffusion or heat equation. 
By this coordinate transformation, the reversible effects of advection are separated
from the irreversible joint effects of advection and diffusion. In this framework, LCSs
express themselves as (boundaries of) metastable sets under the Lagrangian diffusion process. 
In the case of spatially homogeneous isotropic diffusion, averaging the time-dependent
family of Lagrangian diffusion operators yields Froyland's dynamic Laplacian.
In the associated geometric heat equation, the distribution of heat is governed by the
dynamically induced intrinsic geometry on the material manifold,
to which we refer as the \emph{geometry of mixing}.
We study and visualize this geometry in detail, and discuss connections between
geometric features and LCSs viewed as \emph{diffusion barriers} in two numerical examples. Our approach facilitates the
discovery of connections between some prominent methods for coherent structure
detection: the dynamic isoperimetry methodology, the variational geometric approaches
to elliptic LCSs, a class of graph Laplacian-based methods and the effective diffusivity
framework used in physical oceanography.
\end{abstract}



\newpage

\section{Introduction}

Understanding the distribution of physical quantities by advection--diffusion
is of fundamental importance in many scientific disciplines,
including turbulent (geophysical) fluid dynamics and molecular dynamics. Of particular
interest are \emph{coherent structures}, for which there exist many
phenomenological descriptions, visual diagnostics and mathematical
approaches; see \cite{Hadjighasem2017} for a recent review.
In fluid dynamics, coherent structures are often thought
of as rotating islands of particles with regular motion, which move
in an otherwise turbulent background \cite{McWilliams1984,FazleHussain1986,Provenzale1999,Haller2013a}.
In molecular dynamics, coherent structures (or almost-invariant sets) are
thought of as conformations, i.e., sets of configurations of the molecule which
are stable on time scales much larger than those of molecular oscillations \cite{Schutte1999,Schutte2003}.

In the last years, there has been an explosion of coherent structure detection
methods based on flow information. Relying on flow information appears to be a
necessary step in non-auto\-no\-mous/\-un\-stea\-dy velocity fields, since
instantaneous velocity snapshots and their streamlines are no longer
conclusive for material motion, in contrast to the autonomous/steady case.
Nevertheless, the appearance of these methods is very different at
first sight: In the category of variational approaches some methods
require preservation of boundary length \cite{Haller2013a}, minimization of mixing
under the flow \cite{Froyland2010a,Froyland2013} or surface-to-volume ratio
\cite{Froyland2015a,Froyland2015b}. A different class of methods considers
averages of observables along trajectories
\cite{Mezic2010,Budisic2012,Mancho2013,Mundel2014,Haller2016,AlMomani2018} and seeks
coherent structures as sets with similar statistics. Recent (graph) clustering approaches
\cite{Froyland2015,Hadjighasem2016,Banisch2017,Padberg-Gehle2017} assess
coherence based on mutual trajectory distances. Besides, there exist geometric
and topological approaches to coherence; see \cite{Ma2014,Ma2015a} and \cite{Allshouse2012}, respectively.
Comparison studies of these methods have been restricted to simulation case studies
\cite{Allshouse2015,Ma2015,Hadjighasem2017} so far.

Even though many of the above-mentioned approaches focus on different
phenomenological features of coherent structures, often the underlying motivation is
that Lagrangian coherent structures are expected to be material sets which are the least vulnerable to 
(weak) diffusion. By \emph{material} or \emph{Lagrangian sets}, we mean---as usual in continuum
mechanics---flow invariant sets, or, equivalently, fixed sets
of particles. The invulnerability to diffusion is often modeled via some requirement on boundary
deformation under the flow, see, for instance, \cite{Haller2012,Haller2013a,Ma2014,Froyland2015a}.
Despite the intuitive reference to diffusion, all these methods assume a purely advective transport process. In this work, we develop a unifying framework for the
study of coherence from the Lagrangian viewpoint on advection--diffusion and provide new mathematical
connections between Cauchy--Green tensor-based methods developed by Haller and
coworkers \cite{Haller2012,Haller2013a,Farazmand2014a}, the dynamic Laplacian 
methodology by Froyland \cite{Froyland2015a,Froyland2017} and Nakamura's
\emph{effective diffusivity} framework \cite{Nakamura1996,Shuckburgh2003},
adding to the previously found connection between the dynamic Laplacian and the 
probabilistic transfer operator approach \cite{Froyland2015a}; see \Cref{fig:methods} for a schematic overview and more details in \cref{sec:discussion}.

\begin{figure}
\centering
\begin{tikzpicture}[node distance=0.15\textwidth]

\node[draw, text centered, text width =0.17\textwidth] (A) at (0,1.5) {
variational geometric methods};

\node[draw, text centered,text width=0.25\textwidth] (B) at (5.85,0.0) {
dynamic Laplacian, geometric heat flow};

\node[draw, text centered, text width =0.17\textwidth] (C) at (0,-1.5){
graph-Laplace methods};

\node[draw, text centered, text width =0.145\textwidth] (D) at (11.7,1.5){
probabilistic transfer operator};

\node[draw, text centered, text width=0.145\textwidth] (E) at (11.7,-1.5){effective diffusivity};

\draw [<->, color=black, line width=2pt] (A) edge (B) (B) edge (D) (B) edge (C);

\draw [<->, color=black, line width=2pt] (A) -- (B) node [pos=0.55, above] {\rotatebox{-13.5}{\small \Cref{sec:Haller_methods}}};
\draw [<->, color=black, line width=2pt] (C) -- (B) node [pos=0.45, above] {\rotatebox{13.5}{\small \Cref{sec:discretization}}};
\draw [<->, color=black, line width=2pt] (B) -- (D) node [pos=0.5,above] {\rotatebox{13.5}{\small \cite{Froyland2015a,Froyland2017}}};
\draw [<->, color=black, line width=2pt] (B) -- (E) node [pos=0.55, above] {\rotatebox{-13.5}{\small \Cref{sec:gfd}}};
\end{tikzpicture}
\caption{Schematic representation of the connections between different methods
for coherent structure detection.}
\label{fig:methods}
\end{figure}

The intuition of Lagrangian coherence as persistence to diffusion leads us to the
incompressible advection--diffusion equation
(ADE) in Lagrangian coordinates, which is of diffusion-only type; see also
\cite{Press1981,Krol1991,Knobloch1992,Thiffeault2003,Fyrillas2007} for earlier
related approaches. In the Lagrangian frame, we view \emph{Lagrangian coherent 
sets} or \emph{structures (LCSs)} as \emph{metastable sets} under the 
Lagrangian ADE. It turns out that the deformation by advection (in the Eulerian
frame) is equivalent to a deformation of the geometry of the (initial) material manifold, 
i.e., in the Lagrangian frame. This change of perspective from space (Eulerian
frame) to material (Lagrangian frame) solves, by the way, the problem
of separating the reversible effects of pure advection from 
the irreversible effects of advection and diffusion acting together; see \cite{Nakamura1996,Shuckburgh2003}.
Time-averaging of the Lagrangian ADE 
yields an autonomous diffusion-type equation, whose generator coincides---in the case of 
spatially isotropic diffusion---with Froyland's recently introduced \emph{dynamic 
Laplacian} \cite{Froyland2015a}. Froyland's approach is motivated by a dynamic 
analogue to the isoperimetry problem, i.e., the optimal bisection of a manifold, where
optimality is measured with respect to the ratio between the area of the bisection 
surface and the volume of the smaller of the two parts. Our independent and physical
advection--diffusion-based derivation of the self-adjoint dynamic Laplacian establishes a link
to symmetric Markov processes and their metastable decomposition of state
space \cite{Davies1982a,Davies1982,Deuflhard2000,Huisinga2006}.

The Lagrangian averaged diffusion tensor field generates an intrinsic Riemannian
geometry on the material manifold, which we refer to as \emph{geometry of
mixing}; cf.~also \cite{Giona2000} for the use of this terminology,
however, not in an averaging sense. The self-adjoint Laplace operator associated to
the geometry of mixing can be investigated in detail by methods from
semigroup and operator theory \cite{Davies1982,Davies1995}, Riemannian and
spectral geometry \cite{Cheeger1970,Lablee2015} and visualization, i.e., \emph{diffusion tensor imaging (DTI)}.
The technical requirements on the flow are volume-preservation and smoothness, and
on the original material manifold are compactness and boundary regularity that permits 
the formulation of (homogeneous) Neumann boundary conditions. 
Since we are working in a Riemannian geometry setting, we benefit from the
rich intuition about the role of eigenfunctions gathered in applied and computational
harmonic analysis \cite{Coifman2006,Dsilva2015}.

While our theory as presented in the current paper is of continuous type and
theoretically requires arbitrary fine dynamic information, it is strongly related to the
diffusion-maps methodology \cite{Coifman2006}; cf.~also \cite{Banisch2017}. The classic, albeit not exclusive,
application case there is that of \emph{manifold learning}, i.e., the computation of
topological and geometric features (such as intrinsic coordinates) of manifolds
embedded in a Euclidean, usually high-dimensional space. The situation is thus one
of a static manifold. Recently, these ideas have been extended toward including 
dynamics \cite{Giannakis2012,Shnitzer2017,Marshall2018}, taking different approaches than ours.

Another contribution of ours is the following conceptual clarification.
Lagrangian coherent structures (LCSs) are commonly referred to as \emph{transport
barriers}. With its reference to \emph{advection} through the term ``transport'', this
translates then to sets with (near-)zero advective flux. It has been pointed out earlier
\cite{Nakamura1996,Haller2012} that in purely advective flows any material surface
constitutes a transport barrier by flow invariance. Our approach, and its consistency
with many existing LCS methods, clarifies the role of LCSs as \emph{diffusion} or
\emph{mixing barriers}; cf.~\cite{Froyland2013} for an Eulerian analogue.

The paper is organized as follows. We start in \cref{sec:Preliminaries} by recalling
some fundamental concepts from Riemannian geometry, the Laplace operator, its
induced heat flow and metastability. \Cref{sec:Lagrangian_FP}
is devoted to the derivation and discussion of the Lagrangian version
of the advection--diffusion equation (ADE), the definition of Lagrangian coherent
structures in this framework and the derivation of the geometry of mixing.
In \cref{sec:Dynamic-Laplacian}, we study and visualize the geometry of mixing in search
for signatures of diffusion barriers. We close with a discussion of relations to previously developed methods and future directions in
\cref{sec:discussion}. In particular, readers interested in applications in atmospheric 
and oceanic fluid dynamics may find the discussion of connections to the effective 
diffusivity framework in \cref{sec:gfd} of particular interest.

\paragraph*{Notation}

Throughout this paper, we use the following notations.

First, for the symmetric positive-definite matrix representation $G\in\mathbb{R}^{d\times d}$ of a
Riemannian metric $g$ on $M$ (in local coordinates), we denote its
ordered eigenvalues by $0<\mu_{\min}(G)\leq\cdots\leq\mu_{\max}(G),$
and the corresponding eigenvectors (in those coordinates) by $v_{\min}(G),\ldots,$$v_{\max}(G)$.

Second, for any time-dependent map $[0,T]\ni t\mapsto\gamma(t)\in X$,
with $X$ some linear space, we define the time average of $\gamma$
by
\[
\fint_{0}^T\gamma(t)\,\mathrm{d}t\coloneqq\frac{1}{T}\int_{0}^{T}\gamma(t)\,\mathrm{d}t.
\]

\section{Preliminaries}
\label{sec:Preliminaries}

For general references on (weighted) Riemannian manifolds with emphasis on
Laplace operators, heat flows and heat kernels, see \cite{Chavel1984,Rosenberg1997,Grigoryan2009,Jost2011,Lablee2015}.

\subsection{Weighted Manifolds and the Laplace Operator}
\label{sec:weighted_manifolds}

Let $(M,g,\nu)$ be a \emph{weighted manifold}: $M$ a compact, complete, smooth,
connected $d$-dimensional Riemannian manifold, possibly with sufficiently regular boundary
$\partial M$;\footnote{Boundary regularity is generally a delicate matter in PDE
theory. For our purposes, the validity of the \emph{divergence theorem} on $M$ is
most relevant. This theorem is an easy consequence of Stokes's theorem
\cite[Thm.~16.31]{Lee2012}, which can be proven to hold on smooth manifolds with
corners \cite[Thm.~16.25]{Lee2012} in the sense of \cite[p.~415]{Lee2012}.
Throughout, we assume that the manifold $M$ has the so-called \emph{extension
property}, which allows for a rigorous formulation of Neumann boundary conditions,
see \cite[Chap.~7]{Davies1995}, \cite[Sec.~11.5]{Jost2013}.} $g$ a smooth
Riemannian metric (tensor field); and $\nu$ a measure
on $M$ given by integrating indicator functions of measurable sets $A$ with
respect to $\mathrm{d}\nu=\rho\,\mathrm{d}x$,
i.e.,
\[
\nu(A) = \int_M \chi_A(x)\,\mathrm{d}\nu(x) =  \int_A \rho\,\mathrm{d}x.
\]
Here, $\mathrm{d}x$ is the unique (Riemannian) volume form induced by $g$ and
$\rho$ is some smooth positive density \cite{Grigoryan2006}. This gives rise to corresponding
$L^{p}$-spaces over $M$, and in particular to the Hilbert space $L^{2}(M,\nu)$.

For any smooth real-valued function $f\in C^{\infty}(M)$, its \emph{exterior derivative}
$\mathrm{d}f$ is a one-form on $M$, invariantly defined by 
\[
\mathrm{d}f=\frac{\partial f}{\partial x^{i}}\mathrm{d}x^{i}
\]
in local coordinates, where we make use of Einstein's summation convention. Thus, when
viewed as a vector in local coordinates, $\mathrm{d}f$ comprises the partial derivatives of
$f$ in the coordinate directions.\footnote{This is often denoted by $\nabla f$ in the applied literature.}
The \emph{metric tensor} $g_{x}(\cdot,\cdot)$ defines a scalar product on each tangent space $T_{x}M$,
which allows to identify the cotangent space $T_{x}^{*}M$ with the tangent space
$T_{x}M$ via  the \emph{canonical/musical isomorphism}, see
\cite[p.~342]{Lee2012}. In local coordinates, this isomorphism is given by the
inverse Gramian matrix $G^{-1}$, $G$ the matrix representation of the tensor $g$.
In particular, one has
\[
G^{-1}\mathrm{d}f=\grad_gf,
\]
the \emph{gradient of $f$} w.r.t.~the metric $g$. Moreover, the \emph{divergence of a vector field $V$} may be defined implicitly via the Lie derivative of the volume form $\mathrm{d}\nu$ in the direction of $V$,
\[
\mathcal{L}_V(\mathrm{d}\nu)=\divergence_{\nu}(V)\mathrm{d}\nu.
\]
Intuitively, the divergence measures the rate of expansion of volume along the flow
induced by $V$. In particular, if $\divergence(V)=0$, volume is preserved by the flow of $V$.

Finally, the \emph{Laplace operator $\Delta_{g,\nu}$} is defined by
\[
\Delta_{g,\nu} f\coloneqq \divergence_{\nu}\grad_g f = \divergence_{\nu}\,D\,\mathrm{d}f,\qquad D\coloneqq G^{-1}.
\]
Its \emph{weak formulation} takes the form \cite{Grigoryan2006,Grigoryan2009}
\begin{equation}\label{eq:weak_form_weighted}
\left\langle f,\divergence_{\nu}\grad_gh\right\rangle _{0,\nu} = -\int_M g^{-1}(\mathrm{d}f,\mathrm{d}h)\,\mathrm{d}\nu= -\int_M g(\grad_gf,\grad_gh)\,\mathrm{d}\nu,
\end{equation}
where \cref{eq:weak_form_weighted} is known as \emph{Green's formula}. Here, $f,h$ need to satisfy homogeneous Dirichlet or
Neumann boundary conditions \cite{Grigoryan2006}, and $g^{-1}$ denotes the
\emph{dual metric (to $g$)}, which is the pullback of $g$ by the canonical
isomorphism. By construction, the isomorphism is an isometry and, hence, one has
$\lVert \grad_g f(x)\rVert_g = \lVert \mathrm{d}f(x)\rVert_{g^{-1}}$ for any $x\in M$ and any smooth $f$.

It is well known \cite{Grigoryan2006} that $-\Delta_{g,\nu}$ can be represented in local coordinates as
\begin{equation*}
g^{ij}(x)\frac{\partial}{\partial x^{i}}\frac{\partial}{\partial x^{j}}+b^{i}(x)\frac{\partial}{\partial x^{i}}+c(x),
\end{equation*}
where $g^{ij},b^{i},c$ are smooth (real) coefficients and $\left(g^{ij}\right)_{ij}=G^{-1}$ is
symmetric and uniformly positive definite. That is, the principal symbol $g^{ij}(x)\xi_{i}\xi_{j}$ of $-\Delta_{g,\nu}$ satisfies
\[
g^{ij}(x)\xi_{i}\xi_{j}\geq\gamma\lvert\xi\rvert_{g_{x}}^{2},\qquad(x,\xi)\in T^{*}M,
\]
for some $\gamma>0$. In fact, one even has $g^{ij}(x)\xi_{i}\xi_{j}=\lvert\xi\rvert_{g_{x}}^{2}$ by definition. As a consequence, $\Delta_{g,\nu}$ is a \emph{(uniformly) elliptic} second-order differential operator on $M$.

One important property of the Laplace operator is its \emph{invariance under
(volume-preserving) isometries}, see \cite[Sec.~4.2]{Grigoryan2006}: We call
$T\colon N\to M$ an isometry between two weighted manifolds $(N,\tilde{g},\tilde{\nu})$ and
$(M,g,\nu)$ if $T$ is a diffeomorphism and $\tilde{g}=T^*g$ ($\tilde{g}$ is the pullback metric) and $\nu=T_*\tilde{\nu}$ ($\nu$ is the pushforward measure). For such isometries one has
\begin{equation}
\Delta_{\tilde{g},\tilde{\nu}}T^{*}=T^{*}\Delta_{g,\nu}, \qquad\textnormal{or, explicitly}\qquad
\Delta_{\tilde{g},\tilde{\nu}}(f\circ T)= \left(\Delta_{g,\nu}f\right)\circ T,\label{eq:Laplace_trafo}
\end{equation}
for $f\in C^\infty(M)$. This directly implies the coordinate-independence of eigenvalues and eigenfunctions of $\Delta_{g,\nu}$, when interpreting $N$ as a global reparametrization of $M$.

Following a well-established procedure, see \cite[Chap.~3]{Jost2011} for the case
when $\mathrm{d}\nu=\mathrm{d}x$ and \cite[Theorem 2.2]{Grigoryan2006},
$\Delta_{g,\nu}$---defined on smooth functions---can be uniquely extended to a \emph{self-adjoint}
non-positive definite operator on $L^2(M,\nu)$ by the Friedrichs extension. Indeed, 
Green's formula, \cref{eq:weak_form_weighted}, implies that $\Delta_{g,\nu}$ is
$L^2$-symmetric on the domain of classically smooth functions $C^\infty(M)$.
Density of $C^\infty(M)$ in $L^2(M,\nu)$ follows directly from the observation that
$L^2(M,\nu)$ and $L^2(M,\mathrm{d}x)$ are isometric via the unitary transformation
$U_{\sqrt{\rho}}$, $f\mapsto \sqrt{\rho}f$, which leaves $C^\infty(M)$ invariant;
cf.~\cite[Sec.~4.2]{Davies1989}. As a negative semi-definite, self-adjoint elliptic
second-order differential operator, $\Delta_{g,\nu}$ has the following well-known spectral properties.

\begin{prop}[{\cite[Sect.~2.2]{Grigoryan2006}, see also \cite[Thm.~3.2.1]{Jost2011}}]\label{thm:Laplace_spectral}
The operator $\Delta_{g,\nu}$ has
\begin{enumerate}
\item[(i)]  purely discrete, non-positive spectrum $0=\lambda_1\geq\lambda_2\geq\ldots$, and eigenvalues accumulate only at $-\infty$;
\item[(ii)] pairwise $L^2(M,\nu)$-orthogonal eigenspaces;
\item[(iii)] $C^\infty$-smooth eigenfunctions, which form a complete basis of $L^2(M,\nu)$.
\end{enumerate}
The \emph{harmonic functions}, i.e., eigenfunctions corresponding to the 0-eigenvalue, are constant.
Hence, there are as many linearly independent harmonic functions as $M$ has connected components.
\end{prop}

The last means that the multiplicity of the 0-eigenvalue equals the number of connected components
of $M$, if one allowed $M$ to have multiple connected components; we come back to this
issue in \cref{sec:metastability}.

\subsection{Heat Flows}\label{sec:heat-flow}

Given an elliptic, non-positive second-order differential operator $H$
(equipped with zero Neumann boundary condition if $M$ has boundary) such as the
Laplace operator on $M$, the \emph{(infinitesimal) generator},
$\left(\exp(tH)\right) _{t\geq0}$ is an analytic semigroup
of bounded operators defined on $L^{2}(M,\nu)$ \cite{Davies1982a,Davies1995a},
and $u(t)=\exp(tH)u_0$, $u_0\in L^{2}(M)$, is the unique solution of the generalized heat
equation
\begin{equation*}
\frac{\mathrm{d}}{\mathrm{d}t}u(t)=Hu(t),\qquad u(0)=u_0;
\end{equation*}
see, e.g., \cite[Chap.~5, Sec.~1.4]{Davies1989}, and \cite[Sec.~3]{Grigoryan2006}.
The semigroup $(\exp(tH))_{t\geq 0}$ is called the \emph{heat flow} generated by $H$.
By the spectral mapping theorem, we have
\[
\sigma\left(\exp(tH)\right) = \exp\left(t[\sigma(H)]\right),
\]
with corresponding eigenprojections. In other words, it suffices to study the spectrum
and the eigenprojections of the generator $H$ to determine subspaces which are
invariant under the heat flow $\left(\exp(tH)\right)_{t\geq 0}$. For $H=\Delta_{g,\nu}$,
the heat flow maps $L^2$ functions to $C^{\infty}$ functions, the heat kernel 
is symmetric and, hence, the heat flow is a family of self-adjoint operators on
$L^2(M,\nu)$, $\exp(t\Delta_{g,\nu})1_M=1_M$, which are positivity-preserving \cite[Thms.~3.1, 3.3]{Grigoryan2006}.

If we consider a Hölder-continuous curve $t\mapsto H(t)$ of elliptic
second-order differential operators, the unique solution of
\[
\frac{\mathrm{d}}{\mathrm{d}t}u(t)=H(t)u(t),\qquad u(0)=u_0,
\]
is given by $u(t)=U_{H}(t,0)u_0,$ where the \emph{generalized heat process}
$\left\{ U_{H}(t,s)\right\} _{t\geq s\geq0}$ is the non-autonomous
parabolic solution operator generated by $H$, see \cite[Chap.~II]{Amann1995a}.
In particular, it satisfies
\begin{align*}
U_{H}(t,t) & =\id_{L^{2}(M)}, &  & \text{and} & U_{H}(t,\tau)U_{H}(\tau,s) & =U_{H}(t,s), & \text{for } & s\leq\tau\leq t,
\end{align*}
and one has $U_{H}(t,s)=\exp((t-s)H)$ if $H$ is time-translation invariant.
The integral kernel $u_{H}(t,s;\cdot,\cdot)$ of $U_{H}(t,s)$,
\[
U_{H}(t,s)\psi(x)=\int_{M}u_{H}(t,s;x,y)\psi(y)\,\mathrm{d}x,
\]
is called the \emph{heat kernel} of $H$.

\subsection{Metastability and Metastable Decompositions}\label{sec:metastability}

In the literature, it is usually assumed that the manifold under consideration is
connected, as one may otherwise study each connected component individually. In
the spirit of spectral geometry, however, i.e., the study of geometric properties of
Riemannian manifolds by means of spectral properties of the Laplace operator and
its induced heat flow, one important question is the recognition of manifolds that
are connected but ``nearly decomposable''. Such a phenomenon is closely related to the concept
of metastability, which we recall in the following.

First, we recall the seminal work by Davies \cite{Davies1982a,Davies1982} on
metastable states in positivity-preserving contraction semigroups $\left(\exp(tH)\right)_{t\geq 0}$.
The prototype example is given by $H = \Delta$, i.e., a heat flow on a Riemannian manifold,
equipped with homogeneous Neumann boundary condition if $M$ has a non-empty boundary,
as discussed in \cref{sec:heat-flow}. Such operators $H$ have real, non-positive spectrum $0=\lambda_1>\lambda_2\geq \dots$ with a non-degenerate eigenvalue 0, accumulating only at $-\infty$.

For reference, let us formulate some statements:
\begin{enumerate}
\item the first non-trivial eigenvalue $\lambda_2$ of $H$ is much smaller than the second $\lambda_3$;
\item there exists a subset $M_1\subset M$ such that the 
eigenfunction $u_1$ associated with $\lambda_1$ is close (in $L^2$) to some linear combination
of the two indicator functions $1_{M_1}$ and $1_{M\setminus M_1}$;
\item there exists a subset $M_1\subset M$ such that $\lVert\exp(tH)1_{M_1}-1_{M_1}\rVert_1$ is small,
i.e., the indicator function $1_{M_1}$ is slowly evolving away from itself under the heat flow $(\exp(tH))_t$, where
the distance is measured in $L^1$.
\end{enumerate}
Now, \cite[Thm.~7]{Davies1982} and \cite[Thms.~3 \& 5]{Davies1982a} show that,
assuming (1), then (2) and (3) follow, and conversely, assuming (3), then (1) and (2)
follow. Hence, statements (1) and (3) are ``qualitatively equivalent'' \cite[p.~139]{Davies1982a}.
Notably, Davies does not give a direct, explicit definition of metastable sets,
but rather justifies to call $M_1$ a \emph{metastable set} due to the estimates underlying (1)-(3).
Therefore, metastable sets are directly linked to spectral properties of the generator of
the symmetric Markov semigroup. Davies also points out that, having a candidate metastable
set $M_1$, any other set sufficiently close to $M_1$ would admit the same properties, and
hence metastable sets are invariably non-unique.

In \cite{Davies1982}, the above analysis is extended to the physically relevant case, which is of interest in classic LCS applications,
when the generator $H$ has $n$ very small eigenvalues, followed by a significant gap to the next eigenvalue:
\[
0=\lambda_1> \cdots\geq\lambda_n=-\varepsilon, \quad \lambda_{n+1} = O(1).
\]
Similarly to the above, the corresponding first $n$ eigenfunctions are then roughly given as
the linear combination of $n$ indicator functions $1_{M_1},\ldots,1_{M_n}$, where the sets
$M_i$ decompose $M$.

In the context of reversible Markov chains on finite state spaces (e.g., space--time
discretizations of the Markov semigroup above), Davies work has been
adopted and extended as follows.

First, \cite{Deuflhard2000,Deuflhard2005} have studied weak perturbations of
reducible Markov chains. In a geometric heat-flow context, the latter correspond to
multiple connected components, and the weak perturbations introduce weak coupling
between them. In this framework, a perturbation analysis in terms of an
explicit small perturbation parameter $\varepsilon$ on the dominant eigenvalues and
their eigenfunctions is performed and a cluster extraction algorithm called PCCA (or PCCA+)
is devised. In the continuous state space setting as considered by Davies, a corresponding
construction of a perturbation $M_{\varepsilon>0}$ that turns a disconnected reference manifold $M_0$ into a single manifold with ``weakly connected'' components [as done for finite reversible Markov chains in \cite{Deuflhard2000}] appears to be a
very challenging research problem, since both the coupled and the decoupled configurations of 
$M_{\varepsilon\geq 0}$ need to be embedded into a single ambient manifold, on which the Laplace operators and their spectra can be compared.

Another important contribution is \cite{Huisinga2006}, adopting the approach of
\cite{Dellnitz1999} of addressing metastability (called \emph{almost-invariance}
there) of sets via having high internal transition probability (or, equivalently, low exit
probability) as measured by the \emph{(almost) invariance ratio}
\begin{equation*}
p(M_1,M_1) = \langle \mathcal{P}1_{M_1},1_{M_1}\rangle_2/\lVert 1_{M_1}\rVert_2^2.
\end{equation*}
Here, $\mathcal{P}$ is a \emph{Perron--Frobenius/transfer operator}, which can be
thought of as $\exp(tH)$ for some $t>0$ in Davies's framework. In a different line of
research, this approach has been used to detect almost-invariant behavior in
deterministic, finite-dimensional dynamical systems \cite{Froyland2003,Froyland2005}.
The metastability quality of an arbitrary state space decomposition into
$M_1,\ldots,M_n$ is then assessed via $\sum_i p(M_i,M_i)$. Specifically, an upper
bound is given in terms of $\sum_i \lambda_i$, $\lambda_i$ the eigenvalues of
$\mathcal{P}$, and a lower bound in terms of the weighted sum of eigenvalues, where the
weights are determined by the norm of the orthogonal projection of the
eigenfunctions onto the space spanned by the indicator functions $1_{M_i}$. 
Roughly speaking, if the eigenfunctions appear to be almost constant on the $M_i$, 
then the lower bound on metastability is close to the upper bound. In this approach,
there is no a priori assumption on the proximity of dominant eigenvalues to 1. The
provable metastability measure, however, decreases when including eigenfunctions
associated with non-small eigenvalues. For another quantification approach to metastability in terms of exit times, see \cite{Schutte2003}.

It is generally extremely difficult to optimize the almost-invariance objective function
under the constraint of searching among characteristic functions which decompose
state space. For this reason, the optimization problem is relaxed toward
general densities in $L^2$, and the decomposition constraint is relaxed to
$L^2$-orthogonality of densities. In this form, the optimization problem takes the form
of Courant--Fischer min--max type, and is therefore solved by eigenvalues and
eigenfunctions of the (symmetrized) transfer operator. For the extraction of 
metastable/almost-invariant sets, the theory and methods developed in
\cite{Davies1982,Deuflhard2000,Huisinga2006} can then be applied.

In practice, it is common to use heuristic clustering algorithms such as k-means or
fuzzy c-means to extract state-space decompositions from eigenfunctions, besides
the aforementioned \cite{Dellnitz1999,Deuflhard2000,Deuflhard2005}. On the one
hand, there is theoretical justification for using k-means \cite{Lafon2006} by interpreting the
values $\left(w_k(x)\right)_{k=1,\ldots,n}$ of leading eigenfunctions $\left(w_k\right)_{k=1,\ldots,n}$ as a quasi-isometric embedding into some Euclidean space.
In this \emph{feature space}, k-means then effectively optimizes cluster
attribution with respect to the intrinsic distance of the data.
On the other hand, these clustering algorithms do not address the original optimization problem.

\subsection{Laplace Operator, Heat Flow and Local Averaging}
\label{sec:Laplace_averaging}

There exist many tight connections between a Riemannian geometry on a manifold
(as modeled by a Riemannian metric and its induced Laplace operator) on the one
hand, and the heat flow and heat kernel (induced by the Laplace operator) on the
other hand; see, for instance, \cite{Grigoryan2009,Lablee2015}. Here, we want to recall one very
intuitive connection between short-time heat flows and local averaging.

To this end, consider a Riemannian manifold $(M,g)$ without boundary and with
Riemannian measure. Denote the diffusion operator defined by averaging over
$g$–geodesic $\varepsilon$-balls $B_{\varepsilon}^{g}(x)=\lbrace y\in M;\,\dist_{g}(x,y)
\leq\varepsilon\rbrace$ of radius $\varepsilon$ by $T_{\varepsilon}^{g}$, i.e.,
\[
\left(T_{\varepsilon}^{g}u\right)(x)=\frac{\int_{B_{\varepsilon}^g(x)}u\,\mathrm{d}x}{\int_{B_{\varepsilon}^g(x)}\,\mathrm{d}x}=\frac{1}{\Vol_g\left(B_{\varepsilon}^g(x)\right)}\int_{B_{\varepsilon}^g(x)}u\,\mathrm{d}x.
\]
Then, the results from \cite[Thms.~1 and 2]{Lebeau2010} show that
\begin{equation}\label{eq:expansion}
T_{\varepsilon}^{g} =\id_{L^{2}(M,g)}+\tfrac{\varepsilon^{2}}{2(d+2)}\Delta_{g}+O\left(\varepsilon^{4}\right), \qquad \text{for }\varepsilon\to 0,
\end{equation}
(almost) in the norm-resolvent sense; see \cite{Lebeau2010} for technical details. In particular, the dominant eigenvalues and their eigenprojections
of $\varepsilon^2\Delta_{g}$ converge to the eigenvalues and eigenprojections
of $2(d+2)(T_{\varepsilon}^{g}-\id)$ as $\varepsilon\to 0+$, respecting
multiplicity. This strong result can be interpreted in two ways.

\paragraph{Spectral approximation of the short-time heat flow}

As recalled in \cref{sec:metastability}, meta\-stable sets for the heat flow are detected
by eigenfunctions of the heat flow $U_{\varepsilon\Delta}$. Now, the right-hand
side of \cref{eq:expansion} can be read as the second-order operator expansion of
the heat flow $U_{\varepsilon\Delta}\left(\varepsilon/(2d+4)\right)$,
i.e., for short time intervals of length $O(\varepsilon)$. An understanding of this
short-time heat flow is already instructive for the general heat flow, since
\cref{eq:averaged_FP_eq} is
autonomous and the long-time heat flow is nothing but an iteration of the short-time
heat flow. \Cref{eq:expansion}, i.e., the approximation of the short-time heat flow by
the local geodesic averaging operator, then states that it is instructive to look at the shape
distribution of small geodesic neighborhoods to form an intuition on the action of the
short-time heat flow and thereby, possibly, on dominant metastable sets as identified
from Laplace eigenfunctions; cf.~\cref{sec:metastability}. We make use of this
correspondence in a visual exploration of a non-trivial geometry in \cref{sec:tensor-field}.

\paragraph{Approximation of local geodesic averaging by diffusion}

Alternatively, \cref{eq:expansion} may be interpreted from left to right. On the
left-hand side, we have the compact integral smoothing operator $T^g_\varepsilon$.
This operator is expanded in (non-compact) differential smoothing operators. To
zeroth order, i.e., sending $\varepsilon\to 0$, the integral kernel of $T^g$ becomes
the Dirac delta distribution, whose action is given simply by point evaluation, or, on
the operator level, by the identity operator. At the second-order level, local
averaging is represented by the differential diffusion operator $\Delta$.

\section{Advection--Diffusion in Eulerian and Lagrangian Frames}
\label{sec:Lagrangian_FP}

Let $\left(\mathcal{M},\mathfrak{g},\mathrm{d}x\right)$ be a weighted Riemannian $d$-manifold and
$M\subset\mathcal{M}$, the \emph{fluid domain} or \emph{material manifold}, an embedded $d$-dimensional submanifold equipped
with the induced metric, again denoted by $\mathfrak{g}$. We regard $\mathfrak{g}$ and the (volume) measure $\mathrm{d}x$ as universal objects, in the sense that they do not depend on the physical transport and mixing process that we are going to study. In particular, $\mathrm{d}x$ may be the Riemannian measure induced by $\mathfrak{g}$, which is what we assume henceforth. Thus, $\left(\mathcal{M},\mathfrak{g},\mathrm{d}x\right)$ is a Riemannian manifold, and the reader may simply think of the physical space with physical length and volume.

We consider the transport equation/conservation law
for the scalar quantity $\phi$ associated with the (in general non-autonomous) divergence-free vector field $V$ on $\mathcal{M}$:
\begin{equation}
\partial_t\phi + \divergence(\phi V) = 0, \qquad \phi(0,\cdot) = \phi_0.\label{eq:transport_eq}
\end{equation}
Here and throughout, the divergence is the one induced by the physical volume form 
$\mathrm{d}x$.
As is well known, \cref{eq:transport_eq} may be solved for $\phi$ by means of
the \emph{flow map}, i.e., the solution to the ordinary differential equation
\[
\dot{x} = V(t,x),
\]
a smooth one-parametric family of diffeomorphisms $\Phi^{t}$, $t\in[0,T]$
over $M$,
\begin{align*}
\Phi^{t}\colon M & \to\Phi^{t}[M]\subseteqq\mathcal{M}, & t & \in[0,T], & \Phi^{0} & =\id_{M}.
\end{align*}
Indeed, the solution may be represented pointwise by $\phi(t,x) = \phi\left(0,\left(\Phi^t\right)^{-1}(x)\right)$,
or globally $\phi(t,\cdot) = \mathcal{P}^t(\phi(0,\cdot))$, where $\left(\mathcal{P}^t\right)_t$
is the time-dependent family of Perron--Frobenius operators associated with \cref{eq:transport_eq}.

\subsection{Eulerian Advection--Diffusion Equation}\label{sec:Euler_ADE}

As a starting point, consider the spatial evolution of a scalar density $\phi$ as it is
carried by an incompressible fluid and subject to diffusion, the classic Eulerian advection--diffusion equation (ADE)
\begin{equation}
\frac{D\phi}{Dt}=\partial_{t}\phi+\divergence(\phi V) =\varepsilon\divergence\, D\, \mathrm{d}\phi, \qquad \phi(0,\cdot) = \phi_0.\label{eq:euler_ADE}
\end{equation}
Here, $D/Dt$ denotes the \emph{material/substantial/advective time derivative} used in the fluid dynamics literature, and $D$ is the smooth space--time-dependent diffusion tensor field, pointwise symmetric and uniformly positive-definite. The diffusion tensor field is supposed to model only the directional dependence of diffusivity, whereas $\varepsilon>0$ models the diffusion strength and can be interpreted as the inverse of the dimensionless Péclet number,
which quantifies the strength of advection relative to the strength of diffusion.
The problem of LCS detection is typically considered in purely advective flows, which we relax
here to advection-dominated, weakly diffusive flow regimes, i.e., associated with a large Péclet number.

Note that we do not require the spatial metric $\mathfrak{g}$ in the formulation of
\cref{eq:euler_ADE}. It is well known that the above ADE captures also
\emph{anisotropic} diffusion, i.e., diffusion with direction-dependent diffusivity.
Isotropic diffusion corresponds to a diffusion tensor $D$, which is represented by (a
multiple of) the identity tensor in physical units, i.e., in Riemannian normal
coordinates induced by $\mathfrak{g}$.

Nevertheless, we may actually use the diffusion tensor field $D$ to \emph{define} a
diffusion-adapted metric. This is a classic, but little-known procedure, which seems
to go back to Kolmogorov \cite{Kolmogoroff1937}; see also \cite{Masoliver1987,Lara1995}.
It builds on the duality of the Riemannian and the dual metrics, using the
fact that $D$ transforms like a dual metric tensor; see \cref{sec:weighted_manifolds}.
Assume $D$ has components $D^{ij}$ in local coordinates;
then, we define a metric tensor field $g$ by the symmetric, positive-definite matrix field
$G\coloneqq D^{-1}=(D_{ij})_{ij}$. It seems appropriate to refer to (geodesic) length
measured by $g$ as \emph{effective length}: In effective length units, the
$D$-diffusion is isotropic by construction, i.e., in $g$-normal
coordinates $D$ is represented by the identity tensor. This is achieved by
downscaling/upscaling effective length units (relative to physical length units) in directions of
larger/smaller diffusivity, respectively. 
Note that we do not alter the notion of volume as modeled by $\mathrm{d}x$.
By definition, $D$ is the inverse of the metric Gramian $G$ and therefore may be
interpreted as the canonical isomorphism. Thus, $D \mathrm{d}\phi=\grad_g\phi$ 
models the diffusive flux in physical units, and finally \cref{eq:euler_ADE} may be
rewritten in the form
\begin{equation*}
\frac{D\phi}{Dt}=\partial_{t}\phi+\divergence(\phi V) =\varepsilon\divergence\grad_g\phi = \varepsilon\Delta_{g,\mathrm{d}x}\phi, \qquad \phi(0,\cdot) = \phi_0.
\end{equation*}
By the uniform definiteness assumption on the diffusion tensor field $D$, the Laplace
operator $\Delta_{g,\mathrm{d}x}$ is uniformly elliptic. Using local coordinate
representations, it is easy to see that the physical volume measure $\mathrm{d}x$
has a density $\sqrt{\det(D)\det(\mathfrak{G})}$ w.r.t.~the volume measure induced
by $g$, where $\mathfrak{G}$ is the Gramian of $\mathfrak{g}$.

The Eulerian perspective comes with a couple of challenges. First, if one is interested 
in the evolution of material localized in some non-invariant region $M$, one needs to
solve \cref{eq:euler_ADE} on a sufficiently large spatial domain in $\mathcal{M}$ in order
to cover the entire evolution $\Phi^{t}(M)$ of material initialized
in $M$. This can be challenging in applications to open dynamical
systems such as ocean surface flows. Second, coherent sets computed in this
framework---as done in \cite{Denner2016}, cf.~also \cite{Froyland2017b}---are conceptually of
Eulerian, i.e., space--time, kind, and are \emph{not} material by construction.
This lack of materiality is not, as sometimes stated, due to the addition of diffusion in
phase space, as we show by our theory here.
In particular, such Eulerian structures generally have both diffusive and advective fluxes through their
boundary. It is therefore of interest to study weakly diffusive flows from a material perspective, i.e., in Lagrangian coordinates.

\subsection{Lagrangian Advection--Diffusion Equation}

Next, let us take a look at \cref{eq:euler_ADE} from the Lagrangian viewpoint,
cf.~\cite{Press1981,Krol1991,Knobloch1992,Thiffeault2003,Fyrillas2007}. Formally,
this means that we interpret the scalar density as a function of particles by pulling it
back to time $t=0$ through composition with the flow map $\Phi$. This yields a
Lagrangian scalar density $\varphi=\Phi^*\phi = \phi\circ\Phi$. Additionally, we need
to pull back \cref{eq:euler_ADE} to the material manifold, and thus arrive at its Lagrangian form
\begin{align}
\partial_t\varphi &= \varepsilon\divergence\left(D\Phi(t)^{-1}\cdot D\cdot D\Phi(t)^{-\top}\right)\mathrm{d}\varphi, & \varphi(0,\cdot) = \phi_0.\label{eq:lagr_ADE}
\end{align}
Here, the scalar density $\varphi$ is no longer subjected to an advective drift---in the Lagrangian perspective,
we are following trajectories---but is subject to diffusion
generated by the time-dependent family of pullback diffusion tensors. Specifically, denote
by $G(t)^{-1} \coloneqq D\Phi(t)^{-1}\cdot D\cdot D\Phi(t)^{-\top}$ the diffusion tensor in
Lagrangian coordinates. Then, by duality, the matrix field $G(t)$ determines a time-dependent
family of diffusion-adapted pullback metrics $g(t)$ on $M$, and we may rewrite \cref{eq:lagr_ADE} as
\begin{align*}
\partial_t\varphi &= \varepsilon\Delta_{g(t),\mathrm{d}x}\varphi, & \varphi(0,\cdot) = \phi_0,
\end{align*}
where the Laplace operators $\left(\Delta_{g(t),\mathrm{d}x}\right)_t$ are induced by the 
pullback metric $g(t)\coloneqq(\Phi^t)^*g$ and the physical volume. On an abstract level,
this change of notation corresponds exactly to the transformation behavior of Laplace operators,
\cref{eq:Laplace_trafo}. \Cref{eq:lagr_ADE} can thus be viewed as an inhomogeneous, i.e.,
time-dependent, diffusion equation for Lagrangian scalar densities $\varphi$ on $M$.

\begin{rem}[Pullback metrics]\label{rem:pullbackmetric}
The pullback metric $g(t)$ is well known in the theory of kinematics of deforming continua
by the name \emph{(right) Cauchy--Green strain tensor}; see, for instance,
\cite[p.~356]{Abraham1988}. In the typical case when the space $\mathcal{M}$ is
Euclidean and parameterized by the canonical coordinates $x^1,\ldots,x^d$, the 
pullback metric $g(t)$ has the matrix representation $G(t) = \left(D\Phi(t)\right)^\top\cdot D\Phi(t)$,
where $D\Phi(t)$ is the linearized flow map with entries $(\partial_j\Phi^i)_{ij}$.
In general coordinates, one has $G(t)=\left(D\Phi(t)\right)^\top\cdot G\cdot D\Phi(t)$, where
$G$ is the matrix representation of $g$ in local coordinates on $\Phi^t(M)$. 
\end{rem}

\begin{rem}
We stress that—according to our Lagrangian viewpoint—\cref{eq:lagr_ADE}
is an evolution equation on $M$, even if the flow does not keep
$M$ invariant, i.e., $\Phi^{t}(M)\neq M$.
\end{rem}

We are now in the position to define our main object of interest.

\begin{definition}\label{def:LCS}
\emph{Lagrangian coherent sets} $U\subset M$ are material sets that are
metastable under the time-inhomogeneous heat flow \cref{eq:lagr_ADE}, i.e., the
advection--diffusion equation in Lagrangian coordinates. By \emph{Lagrangian
coherent structure} we mean the boundary of a Lagrangian coherent set. Since both
notions define each other unambiguously, we will abbreviate both simply by LCS
and use them mostly synonymously.
\end{definition}

The metric $g(t)$ is different from $g=g(0)$ unless $\Phi^{t}$ is an
isometry, or, in physical terms, unless $\Phi^{t}$ corresponds to
a solid body motion. Therefore, the Lagrangian diffusion is not isotropic with respect to $\mathfrak{g}$, even if the Eulerian diffusion was. This reflects the fact that
the flow deformation may have pushed two
particles apart or together, and thus their material exchange by diffusion
at some later time point is, respectively, less and more likely; see
\Cref{fig:diffusion_types}. These intuitive heuristics have been formalized and exploited in \cite{Thiffeault2003} to reduce the full-dimensional ADE to a one-dimensional ADE along the most contracting direction.

\begin{figure}
\centering
\input{pullback_geometry.tex}
\caption{Schematic visualization of the pullback geometry and the induced diffusion.
The spatial Euclidean geometry (right) is pulled back to the material manifold
from time $t=0.05$ (left) by the flow map $\Phi_0^{0.05}$ for the rotating double gyre, \Cref{exa:transient1}.
A spatial diffusion with variance $\varepsilon=0.1$ (red circle on the right) is pulled back to the red curve on the left, visualizing diffusion in the pullback metric $g(0.05)$ on the material manifold of equal variance.
As can be seen, the red curve reaches further out than material diffusion with same variance in the original metric $g(0)$ (visualized by the green circle) in some directions, while it does not reach as far in others.
This is due to the deformation by the flow.
Note also the duality to the Eulerian deformation perspective presented in
Welander's classic work on two-dimensional turbulence \cite[Fig.~2]{Welander1955}.
}
\label{fig:diffusion_types}
\end{figure}

Since $-\Delta_{g(t),\mathrm{d}x}$ is elliptic for all $t\in[0,T]$ (each is just an isometric
representative of the elliptic Laplace operator $-\Delta_{g,\mathrm{d}x}$ on the flow image),
the solution of \cref{eq:lagr_ADE} in $L^{2}(M)$ is given by the generalized heat flow
$U_{\varepsilon\Delta_{g(t),\mathrm{d}x}}$ associated with 
$\left(\varepsilon\Delta_{g(t),\mathrm{d}x}\right)_t$, cf.~\cref{sec:heat-flow}.

By construction and as a consequence of volume-preservation by $\Phi$, we have the following result.

\begin{prop}\label{lem:time_dep_LB}
For each $t\in[0,T]$, the operator $\Delta_{g(t),\mathrm{d}x}=\divergence_{\mathrm{d}x}\grad_{g(t)}$ is self-adjoint on $L^{2}(M,\mathrm{d}x)$ and admits the spectral properties stated in \Cref{thm:Laplace_spectral}.
\end{prop}

We see that time dependence enters the definition of $\Delta_{g(t),\mathrm{d}x}$
only through the pullback metric/diffusion tensor field, the physical volume
$\mathrm{d}x$ remains unaltered. For this reason, we henceforth omit the measure
in the notation of the Laplace operator.

\begin{rem}
From a geometric viewpoint, one may consider
the heat flow as a tool to study the geometry of manifolds. Then 
\cref{eq:lagr_ADE} may be interpreted as the heat flow on a time-evolving
manifold. This setting has been studied recently from the diffusion-maps point of 
view in \cite{Marshall2018}.
\end{rem}

\subsection{Metastability in Time-Dependent Processes and Its Approximation\label{sec:Time-pert}}

For two reasons, our definition of LCSs given above does not provide a precise mathematical definition,
but rather merely invokes the intuition underlying metastability. First, as recalled in \cref{sec:metastability},
metastability of Markov processes comes with an intrinsic vagueness, and second, as pointed out in \cite[p.~1]{Koltai2016},
``a straightforward definition of a metastable set in the non-stationary, non-equilibrium case may only be given
case-by-case''.

In \cite{Koltai2016}, the authors build on previous work by Froyland \cite{Froyland2013}, and extract
metastable sets (defined as \emph{coherent sets} in the terminology of \cite{Froyland2013}) from features of
the singular vectors of the solution operator associated to the time-dependent process; cf.~also \cite{Marshall2018}.
In our case, this would correspond to the singular vectors of the generalized heat flow $U_{\varepsilon\Delta_{g(t),\mathrm{d}x}}$
introduced earlier. Note, however, that features of singular vectors---like connected regions of uniformly high absolute value---are
generally not expected to be invariant under the generalized heat flow. While this is conceptually unproblematic in Eulerian
approaches to coherence which do not request flow-invariance of the sets of interest,
it is problematic in a Lagrangian approach, where one seeks Lagrangian structures. Those ought to be, by
definition, invariant in Lagrangian coordinates.

To enforce invariance, one may consider extracting LCSs from features of the eigen\-func\-tions of the generalized
heatflow. The issue with this approach then is that the generalized heatflow is a family of non-self-adjoint operators,
whose eigenvalues therefore are not necessarily real. It remains unclear how to interpret complex eigenvalues of modulus
almost 1, and their associated eigenfunctions, in a finite-time context.

Another, simplifying \emph{ad hoc} approach is to approximate the generalized heatflow by an autonomous heat flow, whose generator
is obtained from averaging the time-dependent generators \cite{Press1981,Krol1991,Knobloch1992}, i.e.,
\begin{equation*}
\overline{\Delta}\coloneqq\fint_{0}^{T}\Delta_{g(t)}\,\mathrm{d}t.
\end{equation*}
We defer a rigorous convergence study of the two heatflows in the vanishing-diffusivity limit, $\varepsilon\to 0$, to the
forthcoming \cite{Karrasch2020a}; see \cite{Krol1991} for related work.

By the linearity of the divergence and with the time-average of the pullback diffusion tensors,
\[
\bar{g}^{-1}\coloneqq\fint_0^T g(t)^{-1}\,\mathrm{d}t,
\]
$\overline{\Delta}$ takes the form of a Laplace operator on $M$, i.e.,
\begin{equation*}
\overline{\Delta} = \divergence_{\mathrm{d}x}\grad_{\bar{g}} = \Delta_{\bar{g},\mathrm{d}x}= \Delta_{\bar{g}},
\end{equation*}
where $\bar{g}$ denotes the metric tensor induced by $\bar{g}^{-1}$.
$\Delta_{\bar{g}}$ is then a volume-based diffusion operator again.

Following a different line of reasoning, Froyland introduced the
operator $\overline{\Delta}$ recently in \cite{Froyland2015a} and
coined it \emph{dynamic Laplacian}; cf.~also \cite{Froyland2017}. From our ADE
point of view, the dynamic Laplacian of \cite{Froyland2015a} can be obtained
as the time-average of pullbacks of isotropic Eulerian Laplace operators; cf.~\cref{eq:Laplace_trafo}.

The following proposition holds due to the fact that $\overline{\Delta}$ is the Laplace
operator $\Delta_{\bar{g}}$ associated to the weighted manifold $(M,\bar{g},\mathrm{d}x)$. Notably,
ellipticity follows from uniform bounds on the continuously differentiable flow map 
defined on a compact space-time manifold $M\times [0,T]$.

\begin{prop}[{cf.~\cite[Thm.~4.1]{Froyland2015a}}]
\label{lem:properties_dyn-lapl}
The operator $\overline{\Delta}$ is self-adjoint on $L^{2}(M,\mathrm{d}x)$ and admits the spectral properties stated in \Cref{thm:Laplace_spectral}.
\end{prop}

In the case when $M$ has non-empty boundary, the operator is---as before---equipped
with the natural (w.r.t.~$\bar{g}$) Neumann boundary condition. This can be interpreted as the average of pullbacks
of zero Neumann boundary conditions, see \cite{Froyland2015a,Froyland2017}.

The metric $\bar{g}$ endows the material manifold $M$ with a Riemannian metric,
that encodes---in an averaged sense---the diffusion as it is observed from a Lagrangian perspective. This allows for
(i) a static visualization and exploration and (ii) the application
of many well-established techniques and tools from geometric spectral
analysis and spectral geometry \cite{Jost2011,Lablee2015}, as well as harmonic analysis.

Finally, we propose to approximate LCSs, i.e., metastable sets of the non-autonomous
Lagrangian ADE \eqref{eq:lagr_ADE} by metastable sets for the autonomous
Lagrangian evolution equation
\begin{equation}
\partial_{t}\varphi=\varepsilon\Delta_{\bar{g}}\varphi.\label{eq:averaged_FP_eq}
\end{equation}
By the spectral relation between heat flow and generator, cf.~\cref{sec:heat-flow}, this boils down to a spectral analysis of
the generating Laplace operator $\Delta_{\bar{g}}$. Extracting metastable sets from eigenfunctions of the heat flow/generator
as described in \cref{sec:metastability} then coincides with the procedure put forward in \cite{Froyland2018},
where a couple of example flows are analyzed.

\section{Geometry of Mixing}
\label{sec:Dynamic-Laplacian}

An important outcome of the time averaging of the pullback Laplacians is the
geometric structure on the material manifold $M$ given by the \emph{harmonic mean metric} $\bar{g}$.
This is a material geometry typically different from all material geometries induced by
the configuration of the material in space, i.e., its embedding in space as a submanifold.
The aim of this section is to study this new weighted (Riemannian) manifold as
well as the properties of the induced Laplace operator $\Delta_{\bar{g}}$ and
its heat flow. We will refer to this geometry as the \emph{geometry of mixing}, thereby
reviving and generalizing an earlier related approach \cite{Giona2000}, which
refers to the (single) pullback geometry under one flow map of (typically chaotic) flows.

More specifically, we wish to find signatures of coherent and incoherent
dynamics, and of the boundary between them, in the static geometry of
mixing. We do so by comparing characteristics of the diffusion tensor field
relative to the physical geometry $\left(M,\mathfrak{g}\right)$ and in $\mathfrak{g}$–orthonormal
coordinates $x$. In the Euclidean setting, these are the canonical $x^{i}$-coordinates.

By choosing a reference geometry relative to which we study the deformed
$\bar{g}$–geometry our analysis appears to be somewhat reference dependent.
An analogous approach, however, is common in continuum mechanics, where
deformed configurations are analyzed relative to a reference configuration
\cite{Truesdell2004}. Eventually, the spectrum and the eigenprojections
of $\Delta_{\bar{g}}$—the basis of our coherent structure detection
method—are intrinsic and independent of representations w.r.t.~the reference
configuration. Notably, our geometric construction is observer-independent, or,
equivalently, objective, since Euclidean changes of observer do not change the
notions of length and volume, and the diffusion tensor field $D$ is given intrinsically.

\subsection{Lagrangian Averaged Diffusion Tensor Imaging}\label{sec:tensor-field}

In this section, we explore visually the averaged diffusion tensor field in search of
signatures of Lagrangian coherent structures. The visualization of
second-order tensor fields, and diffusion tensor fields in particular, is referred to as
\emph{diffusion tensor imaging (DTI)} and is a well-established
and active field of research in scientific visualization; see, for instance, \cite{Le-Bihan2001}
for a brief review. We denote the \emph{Lagrangian averaged diffusion tensor} by $\overline{D}$.
Two diffusion phenomena that are of interest in DTI are (i) anisotropy and (ii) some scalar measure of diffusivity.

\subsubsection{Anisotropy and Barriers to Diffusion}\label{sec:anisotropy}

For anisotropy, several scalar measures have been proposed in the diffusion
tensor imaging (DTI) literature, for instance the \emph{volume ratio}, given by
\[
\frac{\prod_i^d\mu_i\left(\overline{D}\right)}{\overline{\mu}\left(\overline{D}\right)^d},
\]
where $\overline{\mu}$ denotes the arithmetic mean of the eigenvalues. This quantifies
the volume of the diffusion ellipsoid relative to the volume of the sphere with radius
$\overline{\mu}$. It takes values between $0$ and $1$, where $1$ corresponds to
isotropy, and $0$ corresponds to a lower-dimensional, degenerate ellipsoid and
hence strong anisotropy. Recall that both measures are not intrinsic quantities of
$\bar{g}$, but are determined by viewing the $\overline{D}$--diffusion tensor in 
$\mathfrak{g}$--orthonormal coordinates.

Generally, for a point $p\in M$ and a $\mathfrak{g}$-unit direction $v$, the $\bar{g}$-norm
of $v$ corresponds to the inverse effective $\overline{D}$-diffusivity in $v$-direction. We are
now looking for a canonical, i.e., diagonalized representation of the diffusion tensor
$\overline{D}$ in physical $\mathfrak{g}$-unit directions. This can be achieved simply 
by computing the eigendecomposition of $\overline{D}$, assuming that the matrix
representation of $\mathfrak{g}$ in the chosen coordinates is the identity. The eigenvalues of
$\overline{D}$ then correspond to the characteristic diffusivities, attained in the
directions of the eigenvectors. In other words, the direction $v_{\max}\left(\overline{D}\right)$
associated with $\mu_{\max}\left(\overline{D}\right)$ corresponds
to the direction of strongest (or fastest) diffusion, and, by duality, to
$v_{\min}\left(\overline{G}\right)$, i.e., the direction which is
most strongly compressed under the change of metric from $\mathfrak{g}$ to
$\bar{g}$, see \Cref{fig:geoballs}. The connection between short-time heat flows and
averaging on small geodesic balls (w.r.t.~the intrinsic geometry, here given by
$\bar{g}$)---recalled in \cref{sec:Laplace_averaging}---indicates that this visualization procedure
may be indicative for the action of the heat flow and the location of metastable states as
identified from spectral information.

\begin{figure}
	\centering
	\begin{tikzpicture}[x=0.5cm, y=0.5cm]
	\draw [color=black, line width=1.5pt] (0,0) circle [radius=6];
	
	\draw[color=red, line width=1.5pt] (0,0)  ellipse [rotate=45, x radius=7, y radius=3];
	\draw [color=red, ->,line width=1.5pt,>=stealth] (0,0) -- (4.95,4.95) node [pos=0.6,above, sloped] {$\scriptstyle\mu_{\max}\left(\overline{G}\right)$};
	\node [color=red] at (6.6,5.4) {$v_{\max}\left(\overline{G}\right)$};
	\node [color=red] at (-3.45,2.75) {$v_{\min}\left(\overline{G}\right)$};
	\draw [->,line width=1.5pt,>=stealth] (0,0) -- (6,0) node [pos=0.42,below,fill=white] {$\scriptstyle\mu_{\max}(G)=\mu_{\min}(G)=1$};
	\draw [color=red, ->,line width=1.5pt,>=stealth] (0,0) -- (-2.12,2.12) node [pos=0.47,below=2pt, sloped] {$\scriptstyle\mu_{\min}\left(\overline{G}\right)$};
	\end{tikzpicture}
	\caption{Schematic visualization of the $g$– and $\bar{g}$–unit circles (black
		and red, resp.) in $\mathfrak{g}$–orthonormal coordinates, cf.\ \Cref{fig:diffusion_types}. Note that $\bar{g}$–diffusion
		is fastest in the direction of $v_{\min}\left(\overline{G}\right)$
		because the $\bar{g}$–distance is the shortest on $g$–circles.
		The corresponding diffusion coefficient in that direction is $\mu_{\max}\left(\overline{D}\right)=1/\mu_{\min}\left(\overline{G}\right)$.
		Note also that $\bar{g}$–unit spheres have typically much smaller volume
		than $\mathfrak{g}$–unit spheres because $\mu_{\min}\left(\overline{G}\right)\leq\mu_{\max}\left(\overline{G}\right)<1$ in large regions of $M$, and consequently $\mathfrak{g}$–unit volumes have $\bar{g}$–volume
		smaller than one there, see \cref{sec:density}.}
	\label{fig:geoballs}
\end{figure}

\subsubsection{(Mean) Diffusivity: Lagrangian Effective Diffusivity and Mixing Regions}

Another quantity of interest in the exploration of diffusion tensor fields is a scalar
measure of total (or, mean) diffusivity. It is common to use the trace of the diffusion tensor in the DTI community, i.e., the sum of eigenvalues. For the mean diffusivity, the trace is additionally normalized, i.e., divided by the dimension. In two-dimensional flows,
the trace as a measure of \emph{absolute} diffusivity (in contrast to the relative
strength as measured by anisotropy) at a point is strongly dominated by the maximal
eigenvalue $\mu_{\max}$. It is natural to interpret regions where the scalar diffusivity
field attains high values as \emph{mixing regions}, in which localized scalar densities
are expected to diffuse very quickly, following preferentially the directions of fastest
diffusion, as discussed in \cref{sec:anisotropy}. In the literature, the method of \emph{trajectory
encounter volume} \cite{Rypina2017}, cf.~also \cite{Padberg-Gehle2017}, and Nakamura's
effective diffusivity methodology, see \cref{sec:gfd}, have been proposed
as a way to compute a Lagrangian diffusivity.

\subsubsection{Density}\label{sec:density}

Throughout, we are working in a weighted manifold setting here, in which the volume
measure $\mathrm{d}x$ generally does not coincide with the volume induced by the
diffusion-adapted metric $g$. It is therefore of interest to study the
deviation of the diffusion-adapted intrinsic volume from the physical volume measure.
This amounts to determining the density of $\mathrm{d}x$ w.r.t.~the Riemannian
volume measure in the geometry of mixing.

However, since intuition is based on the physical notion of volume, we
visualize the inverse problem, i.e., the volume density in the geometry of mixing w.r.t.~$\mathrm{d}x$:
\[
\mathrm{d}\bar{g} =\sqrt{\det\overline{G}}\ \mathrm{d}x=\frac{1}{\sqrt{\det\overline{D}}}\mathrm{d}x.
\]

\subsubsection{Numerical Examples}\label{sec:examples}

For simplicity, we assume that the Eulerian diffusion tensor field in the following
examples is given by the identity tensor in the canonical coordinates. The inclusion of 
a space--time-dependent, anisotropic diffusion tensor field is straightforward, see \Cref{rem:pullbackmetric}.
We are going to visualize the geometry of mixing for two commonly studied flow examples.
For an LCS analysis based on spectral data of the dynamic Laplacian for these examples, see \cite{Froyland2018}.

\begin{example}[Rotating double gyre flow \cite{Mosovsky2011}]
	\label{exa:transient1}
	We consider the transient double-gyre flow on the unit square $[0,1]\times[0,1]$, as introduced in \cite{Mosovsky2011}.
	It is given by a time-dependent stream function $\Psi(t,x,y)=(1-s(t))\sin(2\pi x)\sin(\pi y)+s(t)\sin(\pi x)\sin(2\pi y)$,
	$s(t)=t^{2}(3-2t)$, defining the velocity field via
	\begin{equation*}
	\dot{x} =-\frac{\partial\Psi}{\partial y}, \qquad \dot{y} =\frac{\partial\Psi}{\partial x}.
	\end{equation*}
	The integration time interval is $[0,1]$ and the computational grid is 500$\times$500 points. The flow is designed to
	interpolate in time an instantaneously horizontal (at $t=0$) and
	an instantaneously vertical (at $t=1$) double-gyre vector field.
	For our metric computations, we average over 21 pullback metrics from
	equidistant time instances with time step $0.05$. The LCSs as computed
	from a clustering of the dominant eigenfunctions of the dynamic Laplacian $\overline{\Delta}$
	are shown in \Cref{fig:rdg_clustering}; see \cite{Froyland2018} for details.
	For a visual proof of coherence, we provide an advection movie showing the
	evolution of the Lagrangian coherent structures as Supplementary Material 1. 
	
	\begin{figure}
		\centering
		\input{rdg_final.tex}
		\caption{Rotating double-gyre flow: the three-clustering obtained from the second and third eigenfunction
			of $\overline{\Delta}$ at initial (left) and final time (right).}
		\label{fig:rdg_clustering}
	\end{figure}
\end{example}

\begin{example}[Bickley jet flow \cite{Rypina2007}]
	\label{exa:bickley1}
	We consider the Bickley jet flow, as introduced in \cite{Rypina2007}, which is determined by the stream
	function $\psi(t,x,y) =\psi_{0}(y)+\psi_{1}(t,x,y)$, where
	\begin{align*}
	\psi_{0}(y) & =-U_{0}L_{0}\tanh\left(y/L_{0}\right), & 
	\psi_{1}(t,x,y) & =U_{0}L_{0}\sech^{2}(y/L_{0})\Re\left(\sum_{n=1}^{3}f_{n}(t)\exp\left(ik_{n}x\right)\right).
	\end{align*}
	with functions and parameters as in \cite{Rypina2007,Hadjighasem2016}: $f_{n}(t)=\epsilon_{n}\exp\left(-ik_{n}c_{n}t\right)$,
	$U_{0}=62.66 \textrm{m\,s}^{-1}$, $L_{0}=1770 \mathrm{km}$, $k_{n}=2n/r_{0}$, $r_{0}=6.371 \mathrm{km}$,
	$c_{1}=0.1446U_{0}$, $c_{2}=0.205U_{0}$, $c_{3}=0.461U_{0}$, $\epsilon_{1}=0.0075$,
	$\epsilon_{2}=0.15$, $\epsilon_{3}=0.3$; $x$ and $y$ have units of 1000 km and $t$ has unit s.
	The integration time interval is $[0,40]$ days and the computational grid is 800$\times$240 points. 
	We approximate $\bar{g}$ by 81 pullback metrics from equidistant time instances
	with time step $0.5$ days. The LCSs as computed from a clustering of the dominant
	eigenfunctions of the dynamic Laplacian $\overline{\Delta}$ are shown in \Cref{fig:bickley_final};
	see \cite{Froyland2018} again for details.
	As for the previous example, we provide an advection movie showing the
	evolution of the detected sets as Supplementary Material 2 for a visual confirmation of the coherent motion.
	
	\begin{figure}
	\centering
	\begin{tikzpicture}
	
	\begin{axis}[%
	width=0.8\textwidth,
	height=0.24\textwidth,
	at={(1.011in,2.112in)},
	scale only axis,
	axis on top,
	xmin=-0.0125786163522013,
	xmax=20.0125786163522,
	xtick={0,5,10,15,20},
	ymin=-3.01255230125523,
	ymax=3.01255230125523,
	ytick={-3,0,3},
	axis background/.style={fill=white}
	]
	\addplot [forget plot] graphics [xmin=-0.0125786163522013,xmax=20.0125786163522,ymin=-3.01255230125523,ymax=3.01255230125523] {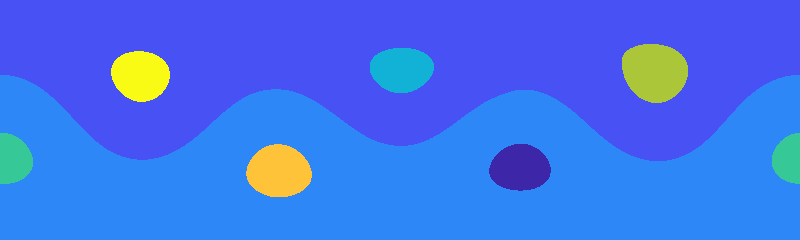};
	\end{axis}
	\end{tikzpicture}%
	\caption{Bickley jet flow: the eight-clustering obtained from the second to eighth eigenfunctions of $\overline{\Delta}$.}
	\label{fig:bickley_final}
	\end{figure}	
\end{example}

\begin{figure}
\centering
\subfloat[Lagrangian averaged diffusion anisotropy]{\label{fig:rdg_anisotropy}
\begin{tikzpicture}
\begin{axis}[%
width=.3\textwidth,
height=.3\textwidth,
at={(1.267in,0.642in)},
xtick={0,0.5,1},
ytick={0,0.5,1},
scale only axis,
point meta min=0,
point meta max=1,
ticklabel style = {font=\small},
every axis/.append style={font=\small},
axis on top,
xmin=0,
xmax=1,
xlabel={$x$},
ymin=0,
ymax=1,
ylabel={$y$},
axis background/.style={fill=white},
colormap={mymap}{[1pt] rgb(0pt)=(0.2422,0.1504,0.6603); rgb(1pt)=(0.25039,0.164995,0.707614); rgb(2pt)=(0.257771,0.181781,0.751138); rgb(3pt)=(0.264729,0.197757,0.795214); rgb(4pt)=(0.270648,0.214676,0.836371); rgb(5pt)=(0.275114,0.234238,0.870986); rgb(6pt)=(0.2783,0.255871,0.899071); rgb(7pt)=(0.280333,0.278233,0.9221); rgb(8pt)=(0.281338,0.300595,0.941376); rgb(9pt)=(0.281014,0.322757,0.957886); rgb(10pt)=(0.279467,0.344671,0.971676); rgb(11pt)=(0.275971,0.366681,0.982905); rgb(12pt)=(0.269914,0.3892,0.9906); rgb(13pt)=(0.260243,0.412329,0.995157); rgb(14pt)=(0.244033,0.435833,0.998833); rgb(15pt)=(0.220643,0.460257,0.997286); rgb(16pt)=(0.196333,0.484719,0.989152); rgb(17pt)=(0.183405,0.507371,0.979795); rgb(18pt)=(0.178643,0.528857,0.968157); rgb(19pt)=(0.176438,0.549905,0.952019); rgb(20pt)=(0.168743,0.570262,0.935871); rgb(21pt)=(0.154,0.5902,0.9218); rgb(22pt)=(0.146029,0.609119,0.907857); rgb(23pt)=(0.138024,0.627629,0.89729); rgb(24pt)=(0.124814,0.645929,0.888343); rgb(25pt)=(0.111252,0.6635,0.876314); rgb(26pt)=(0.0952095,0.679829,0.859781); rgb(27pt)=(0.0688714,0.694771,0.839357); rgb(28pt)=(0.0296667,0.708167,0.816333); rgb(29pt)=(0.00357143,0.720267,0.7917); rgb(30pt)=(0.00665714,0.731214,0.766014); rgb(31pt)=(0.0433286,0.741095,0.73941); rgb(32pt)=(0.0963952,0.75,0.712038); rgb(33pt)=(0.140771,0.7584,0.684157); rgb(34pt)=(0.1717,0.766962,0.655443); rgb(35pt)=(0.193767,0.775767,0.6251); rgb(36pt)=(0.216086,0.7843,0.5923); rgb(37pt)=(0.246957,0.791795,0.556743); rgb(38pt)=(0.290614,0.79729,0.518829); rgb(39pt)=(0.340643,0.8008,0.478857); rgb(40pt)=(0.3909,0.802871,0.435448); rgb(41pt)=(0.445629,0.802419,0.390919); rgb(42pt)=(0.5044,0.7993,0.348); rgb(43pt)=(0.561562,0.794233,0.304481); rgb(44pt)=(0.617395,0.787619,0.261238); rgb(45pt)=(0.671986,0.779271,0.2227); rgb(46pt)=(0.7242,0.769843,0.191029); rgb(47pt)=(0.773833,0.759805,0.16461); rgb(48pt)=(0.820314,0.749814,0.153529); rgb(49pt)=(0.863433,0.7406,0.159633); rgb(50pt)=(0.903543,0.733029,0.177414); rgb(51pt)=(0.939257,0.728786,0.209957); rgb(52pt)=(0.972757,0.729771,0.239443); rgb(53pt)=(0.995648,0.743371,0.237148); rgb(54pt)=(0.996986,0.765857,0.219943); rgb(55pt)=(0.995205,0.789252,0.202762); rgb(56pt)=(0.9892,0.813567,0.188533); rgb(57pt)=(0.978629,0.838629,0.176557); rgb(58pt)=(0.967648,0.8639,0.16429); rgb(59pt)=(0.96101,0.889019,0.153676); rgb(60pt)=(0.959671,0.913457,0.142257); rgb(61pt)=(0.962795,0.937338,0.12651); rgb(62pt)=(0.969114,0.960629,0.106362); rgb(63pt)=(0.9769,0.9839,0.0805)},
colorbar
]
\addplot [forget plot] graphics [xmin=-0.0008998, xmax=1.0008998, ymin=-0.0008998, ymax=1.0008998] {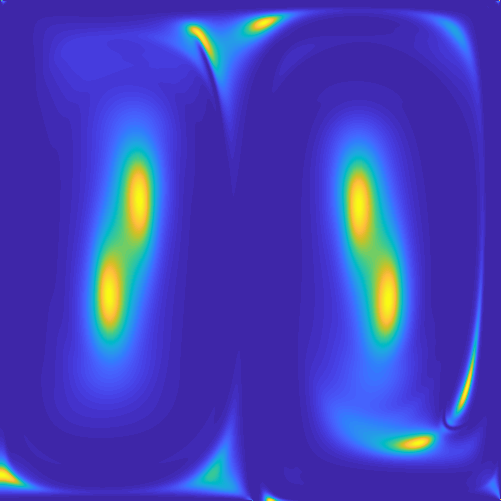};
\end{axis}
\end{tikzpicture}}
\subfloat[Logarithm of the Lagrangian effective diffusivity]{\label{fig:rdg_diffusivity}
\begin{tikzpicture}
\begin{axis}[%
width=.3\textwidth,
height=.3\textwidth,
at={(1.267in,0.642in)},
xtick={0,0.5,1},
ytick={0,0.5,1},
scale only axis,
point meta min=0,
point meta max=7.25467832205814,
axis on top,
ticklabel style={font=\small},
every axis/.append style={font=\small},
xmin=0,
xmax=1,
xlabel={$x$},
ymin=0,
ymax=1,
ylabel={$y$},
axis background/.style={fill=white},
colormap={mymap}{[1pt] rgb(0pt)=(0.2081,0.1663,0.5292); rgb(1pt)=(0.211624,0.189781,0.577676); rgb(2pt)=(0.212252,0.213771,0.626971); rgb(3pt)=(0.2081,0.2386,0.677086); rgb(4pt)=(0.195905,0.264457,0.7279); rgb(5pt)=(0.170729,0.291938,0.779248); rgb(6pt)=(0.125271,0.324243,0.830271); rgb(7pt)=(0.0591333,0.359833,0.868333); rgb(8pt)=(0.0116952,0.38751,0.881957); rgb(9pt)=(0.00595714,0.408614,0.882843); rgb(10pt)=(0.0165143,0.4266,0.878633); rgb(11pt)=(0.0328524,0.443043,0.871957); rgb(12pt)=(0.0498143,0.458571,0.864057); rgb(13pt)=(0.0629333,0.47369,0.855438); rgb(14pt)=(0.0722667,0.488667,0.8467); rgb(15pt)=(0.0779429,0.503986,0.838371); rgb(16pt)=(0.0793476,0.520024,0.831181); rgb(17pt)=(0.0749429,0.537543,0.826271); rgb(18pt)=(0.0640571,0.556986,0.823957); rgb(19pt)=(0.0487714,0.577224,0.822829); rgb(20pt)=(0.0343429,0.596581,0.819852); rgb(21pt)=(0.0265,0.6137,0.8135); rgb(22pt)=(0.0238905,0.628662,0.803762); rgb(23pt)=(0.0230905,0.641786,0.791267); rgb(24pt)=(0.0227714,0.653486,0.776757); rgb(25pt)=(0.0266619,0.664195,0.760719); rgb(26pt)=(0.0383714,0.674271,0.743552); rgb(27pt)=(0.0589714,0.683757,0.725386); rgb(28pt)=(0.0843,0.692833,0.706167); rgb(29pt)=(0.113295,0.7015,0.685857); rgb(30pt)=(0.145271,0.709757,0.664629); rgb(31pt)=(0.180133,0.717657,0.642433); rgb(32pt)=(0.217829,0.725043,0.619262); rgb(33pt)=(0.258643,0.731714,0.595429); rgb(34pt)=(0.302171,0.737605,0.571186); rgb(35pt)=(0.348167,0.742433,0.547267); rgb(36pt)=(0.395257,0.7459,0.524443); rgb(37pt)=(0.44201,0.748081,0.503314); rgb(38pt)=(0.487124,0.749062,0.483976); rgb(39pt)=(0.530029,0.749114,0.466114); rgb(40pt)=(0.570857,0.748519,0.44939); rgb(41pt)=(0.609852,0.747314,0.433686); rgb(42pt)=(0.6473,0.7456,0.4188); rgb(43pt)=(0.683419,0.743476,0.404433); rgb(44pt)=(0.71841,0.741133,0.390476); rgb(45pt)=(0.752486,0.7384,0.376814); rgb(46pt)=(0.785843,0.735567,0.363271); rgb(47pt)=(0.818505,0.732733,0.34979); rgb(48pt)=(0.850657,0.7299,0.336029); rgb(49pt)=(0.882433,0.727433,0.3217); rgb(50pt)=(0.913933,0.725786,0.306276); rgb(51pt)=(0.944957,0.726114,0.288643); rgb(52pt)=(0.973895,0.731395,0.266648); rgb(53pt)=(0.993771,0.745457,0.240348); rgb(54pt)=(0.999043,0.765314,0.216414); rgb(55pt)=(0.995533,0.786057,0.196652); rgb(56pt)=(0.988,0.8066,0.179367); rgb(57pt)=(0.978857,0.827143,0.163314); rgb(58pt)=(0.9697,0.848138,0.147452); rgb(59pt)=(0.962586,0.870514,0.1309); rgb(60pt)=(0.958871,0.8949,0.113243); rgb(61pt)=(0.959824,0.921833,0.0948381); rgb(62pt)=(0.9661,0.951443,0.0755333); rgb(63pt)=(0.9763,0.9831,0.0538)},
colorbar,
]
\addplot [forget plot] graphics [xmin=-0.000977471624266145,xmax=1.00097747162427,ymin=-0.000977471624266145,ymax=1.00097747162427] {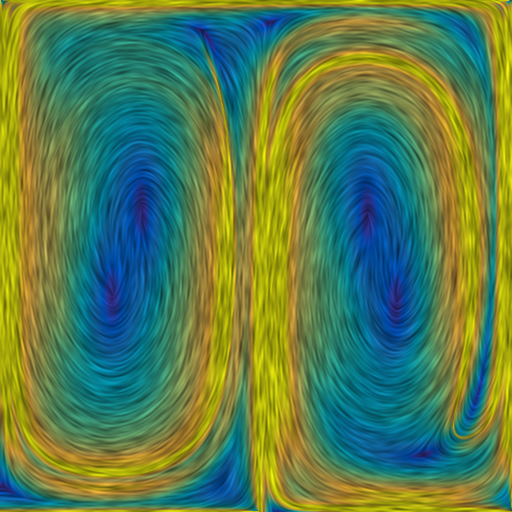};
\end{axis}
\end{tikzpicture}}\\
\subfloat[Volume density]{\label{fig:rdg_density}
\begin{tikzpicture}

\begin{axis}[%
width=.3\textwidth,
height=.3\textwidth,
at={(1.267in,0.642in)},
xtick={0,0.5,1},
ytick={0,0.5,1},
scale only axis,
point meta min=1.32895474792832e-07,
point meta max=0.912257951750063,
axis on top,
xmin=0,
xmax=1,
xlabel={$x$},
ymin=0,
ymax=1,
ylabel={$y$},
axis background/.style={fill=white},
colormap/jet,
colorbar,
]
\addplot [forget plot] graphics [xmin=-0.000977471624266145,xmax=1.00097747162427,ymin=-0.000977471624266145,ymax=1.00097747162427] {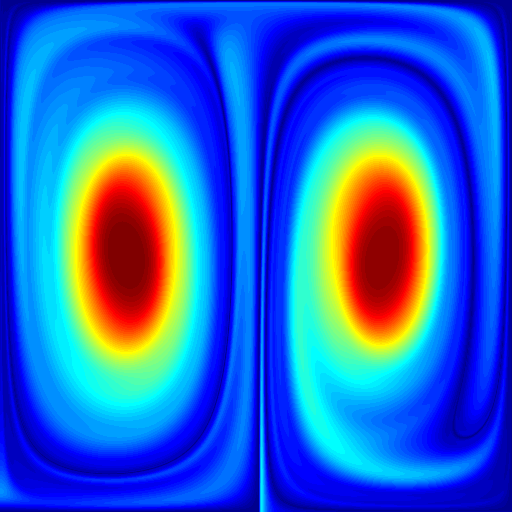};
\end{axis}
\end{tikzpicture}}
\caption{Averaged Lagrangian diffusion tensor field imaging for the transient double gyre flow. (a) Volume ratio as a measure of anisotropy. High values correspond to isotropic diffusion. (b) The texture corresponds to integral curves of the dominant diffusion direction field; the scalar field corresponds to the logarithmic trace of the diffusion tensor, the Lagrangian effective diffusivity. (c) Density of the diffusion-induced volume relative to physical volume.}
\label{fig:gyre_DTI}
\end{figure}

\begin{figure}
\centering
\subfloat[Lagrangian averaged diffusion anisotropy]{\label{fig:Bickley_anisotropy}
\begin{tikzpicture}

\begin{axis}[%
width=.75\textwidth,
height=0.225\textwidth,
at={(2.499in,2.969in)},
scale only axis,
point meta min=0,
point meta max=1,
axis on top,
xmin=0,
xmax=20,
xtick={0,5,10,15,20},
ticklabel style = {font=\small},
every axis/.append style={font=\small},
xlabel={$x$ [Mm]},
ymin=-3,
ymax=3,
ylabel={$y$ [Mm]},
ytick={-3,0,3},
axis background/.style={fill=white},
colormap={mymap}{[1pt] rgb(0pt)=(0.2422,0.1504,0.6603); rgb(1pt)=(0.25039,0.164995,0.707614); rgb(2pt)=(0.257771,0.181781,0.751138); rgb(3pt)=(0.264729,0.197757,0.795214); rgb(4pt)=(0.270648,0.214676,0.836371); rgb(5pt)=(0.275114,0.234238,0.870986); rgb(6pt)=(0.2783,0.255871,0.899071); rgb(7pt)=(0.280333,0.278233,0.9221); rgb(8pt)=(0.281338,0.300595,0.941376); rgb(9pt)=(0.281014,0.322757,0.957886); rgb(10pt)=(0.279467,0.344671,0.971676); rgb(11pt)=(0.275971,0.366681,0.982905); rgb(12pt)=(0.269914,0.3892,0.9906); rgb(13pt)=(0.260243,0.412329,0.995157); rgb(14pt)=(0.244033,0.435833,0.998833); rgb(15pt)=(0.220643,0.460257,0.997286); rgb(16pt)=(0.196333,0.484719,0.989152); rgb(17pt)=(0.183405,0.507371,0.979795); rgb(18pt)=(0.178643,0.528857,0.968157); rgb(19pt)=(0.176438,0.549905,0.952019); rgb(20pt)=(0.168743,0.570262,0.935871); rgb(21pt)=(0.154,0.5902,0.9218); rgb(22pt)=(0.146029,0.609119,0.907857); rgb(23pt)=(0.138024,0.627629,0.89729); rgb(24pt)=(0.124814,0.645929,0.888343); rgb(25pt)=(0.111252,0.6635,0.876314); rgb(26pt)=(0.0952095,0.679829,0.859781); rgb(27pt)=(0.0688714,0.694771,0.839357); rgb(28pt)=(0.0296667,0.708167,0.816333); rgb(29pt)=(0.00357143,0.720267,0.7917); rgb(30pt)=(0.00665714,0.731214,0.766014); rgb(31pt)=(0.0433286,0.741095,0.73941); rgb(32pt)=(0.0963952,0.75,0.712038); rgb(33pt)=(0.140771,0.7584,0.684157); rgb(34pt)=(0.1717,0.766962,0.655443); rgb(35pt)=(0.193767,0.775767,0.6251); rgb(36pt)=(0.216086,0.7843,0.5923); rgb(37pt)=(0.246957,0.791795,0.556743); rgb(38pt)=(0.290614,0.79729,0.518829); rgb(39pt)=(0.340643,0.8008,0.478857); rgb(40pt)=(0.3909,0.802871,0.435448); rgb(41pt)=(0.445629,0.802419,0.390919); rgb(42pt)=(0.5044,0.7993,0.348); rgb(43pt)=(0.561562,0.794233,0.304481); rgb(44pt)=(0.617395,0.787619,0.261238); rgb(45pt)=(0.671986,0.779271,0.2227); rgb(46pt)=(0.7242,0.769843,0.191029); rgb(47pt)=(0.773833,0.759805,0.16461); rgb(48pt)=(0.820314,0.749814,0.153529); rgb(49pt)=(0.863433,0.7406,0.159633); rgb(50pt)=(0.903543,0.733029,0.177414); rgb(51pt)=(0.939257,0.728786,0.209957); rgb(52pt)=(0.972757,0.729771,0.239443); rgb(53pt)=(0.995648,0.743371,0.237148); rgb(54pt)=(0.996986,0.765857,0.219943); rgb(55pt)=(0.995205,0.789252,0.202762); rgb(56pt)=(0.9892,0.813567,0.188533); rgb(57pt)=(0.978629,0.838629,0.176557); rgb(58pt)=(0.967648,0.8639,0.16429); rgb(59pt)=(0.96101,0.889019,0.153676); rgb(60pt)=(0.959671,0.913457,0.142257); rgb(61pt)=(0.962795,0.937338,0.12651); rgb(62pt)=(0.969114,0.960629,0.106362); rgb(63pt)=(0.9769,0.9839,0.0805)},
colorbar
]
\addplot [forget plot] graphics [xmin=-0.0125156445556946, xmax=20.0125156445557, ymin=-3.01255230125523, ymax=3.01255230125523] {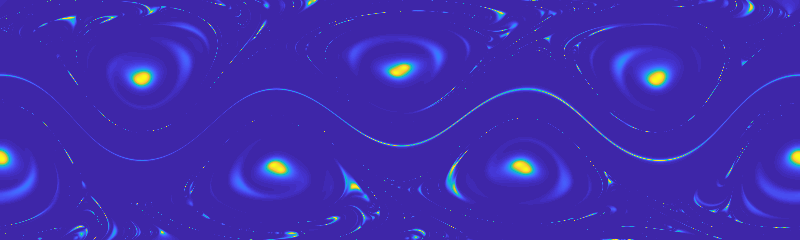};
\end{axis}
\end{tikzpicture}}

\subfloat[Logarithm of the Lagrangian effective diffusivity]{\label{fig:Bickley_diffusivity}
\begin{tikzpicture}

\begin{axis}[%
width=.75\textwidth,
height=0.225\textwidth,
at={(2.499in,2.969in)},
scale only axis,
point meta min=0,
point meta max=13.6391610790367,
axis on top,
xmin=0,
xmax=20,
ticklabel style = {font=\small},
every axis/.append style={font=\small},
xtick={0,5,10,15,20},
xlabel={$x$ [Mm]},
ymin=-3,
ymax=3,
ylabel={$y$ [Mm]},
ytick={-3,0,3},
axis background/.style={fill=white},
colormap={mymap}{[1pt] rgb(0pt)=(0.2081,0.1663,0.5292); rgb(1pt)=(0.211624,0.189781,0.577676); rgb(2pt)=(0.212252,0.213771,0.626971); rgb(3pt)=(0.2081,0.2386,0.677086); rgb(4pt)=(0.195905,0.264457,0.7279); rgb(5pt)=(0.170729,0.291938,0.779248); rgb(6pt)=(0.125271,0.324243,0.830271); rgb(7pt)=(0.0591333,0.359833,0.868333); rgb(8pt)=(0.0116952,0.38751,0.881957); rgb(9pt)=(0.00595714,0.408614,0.882843); rgb(10pt)=(0.0165143,0.4266,0.878633); rgb(11pt)=(0.0328524,0.443043,0.871957); rgb(12pt)=(0.0498143,0.458571,0.864057); rgb(13pt)=(0.0629333,0.47369,0.855438); rgb(14pt)=(0.0722667,0.488667,0.8467); rgb(15pt)=(0.0779429,0.503986,0.838371); rgb(16pt)=(0.0793476,0.520024,0.831181); rgb(17pt)=(0.0749429,0.537543,0.826271); rgb(18pt)=(0.0640571,0.556986,0.823957); rgb(19pt)=(0.0487714,0.577224,0.822829); rgb(20pt)=(0.0343429,0.596581,0.819852); rgb(21pt)=(0.0265,0.6137,0.8135); rgb(22pt)=(0.0238905,0.628662,0.803762); rgb(23pt)=(0.0230905,0.641786,0.791267); rgb(24pt)=(0.0227714,0.653486,0.776757); rgb(25pt)=(0.0266619,0.664195,0.760719); rgb(26pt)=(0.0383714,0.674271,0.743552); rgb(27pt)=(0.0589714,0.683757,0.725386); rgb(28pt)=(0.0843,0.692833,0.706167); rgb(29pt)=(0.113295,0.7015,0.685857); rgb(30pt)=(0.145271,0.709757,0.664629); rgb(31pt)=(0.180133,0.717657,0.642433); rgb(32pt)=(0.217829,0.725043,0.619262); rgb(33pt)=(0.258643,0.731714,0.595429); rgb(34pt)=(0.302171,0.737605,0.571186); rgb(35pt)=(0.348167,0.742433,0.547267); rgb(36pt)=(0.395257,0.7459,0.524443); rgb(37pt)=(0.44201,0.748081,0.503314); rgb(38pt)=(0.487124,0.749062,0.483976); rgb(39pt)=(0.530029,0.749114,0.466114); rgb(40pt)=(0.570857,0.748519,0.44939); rgb(41pt)=(0.609852,0.747314,0.433686); rgb(42pt)=(0.6473,0.7456,0.4188); rgb(43pt)=(0.683419,0.743476,0.404433); rgb(44pt)=(0.71841,0.741133,0.390476); rgb(45pt)=(0.752486,0.7384,0.376814); rgb(46pt)=(0.785843,0.735567,0.363271); rgb(47pt)=(0.818505,0.732733,0.34979); rgb(48pt)=(0.850657,0.7299,0.336029); rgb(49pt)=(0.882433,0.727433,0.3217); rgb(50pt)=(0.913933,0.725786,0.306276); rgb(51pt)=(0.944957,0.726114,0.288643); rgb(52pt)=(0.973895,0.731395,0.266648); rgb(53pt)=(0.993771,0.745457,0.240348); rgb(54pt)=(0.999043,0.765314,0.216414); rgb(55pt)=(0.995533,0.786057,0.196652); rgb(56pt)=(0.988,0.8066,0.179367); rgb(57pt)=(0.978857,0.827143,0.163314); rgb(58pt)=(0.9697,0.848138,0.147452); rgb(59pt)=(0.962586,0.870514,0.1309); rgb(60pt)=(0.958871,0.8949,0.113243); rgb(61pt)=(0.959824,0.921833,0.0948381); rgb(62pt)=(0.9661,0.951443,0.0755333); rgb(63pt)=(0.9763,0.9831,0.0538)},
colorbar,
]
\addplot [forget plot] graphics [xmin=-0.0125156445556946,xmax=20.0125156445557,ymin=-3.01255230125523,ymax=3.01255230125523] {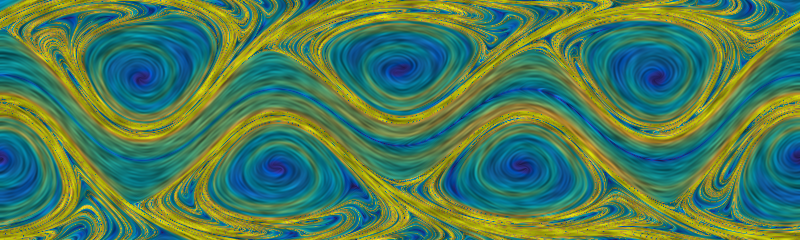};
\end{axis}
\end{tikzpicture}}

\subfloat[Volume density]{\label{fig:Bickley_density}
\begin{tikzpicture}

\begin{axis}[%
width=.75\textwidth,
height=0.225\textwidth,
at={(2.499in,2.969in)},
scale only axis,
point meta min=8.02072779714928e-08,
point meta max=0.994821701535397,
axis on top,
xmin=0,
xmax=20,
xtick={0,5,10,15,20},
xlabel={$x$ [Mm]},
ymin=-3,
ymax=3,
ylabel={$y$ [Mm]},
ytick={-3,0,3},
ticklabel style = {font=\small},
every axis/.append style={font=\small},
axis background/.style={fill=white},
colormap/jet,
colorbar
]
\addplot [forget plot] graphics [xmin=-0.0125156445556946,xmax=20.0125156445557,ymin=-3.01255230125523,ymax=3.01255230125523] {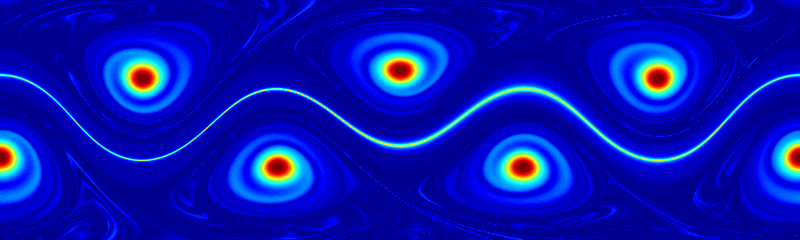};
\end{axis}
\end{tikzpicture}%
}
\caption{Averaged Lagrangian diffusion tensor field imaging for the Bickley jet flow;
see \Cref{fig:gyre_DTI} for the interpretation.}
\label{fig:bickley_DTI}
\end{figure}

Next, we visualize the Lagrangian averaged diffusion tensor field for the two example
flows described in \Cref{exa:transient1,exa:bickley1}. In \Cref{fig:rdg_anisotropy,fig:Bickley_anisotropy}
we show the respective volume ratio fields. In \Cref{fig:rdg_diffusivity,fig:Bickley_diffusivity},
the scalar field shown is the decimal logarithm of the mean diffusivity and it is overlaid by a grayscale
texture whose features are aligned with the dominant diffusion direction field
$v_{\max}\left(\overline{D}\right)$. In the orthogonal direction,
diffusion is weaker typically by several orders of magnitudes.
Finally, in \Cref{fig:rdg_density,fig:Bickley_density} we plot the respective densities $\left(\det\overline{D}\right)^{-1/2}$.

In both \Cref{fig:gyre_DTI,fig:bickley_DTI}, two phenomena are clearly visible.
First, the $\bar{g}$—diffusion gets closest to isotropic diffusion (volume ratio
values close to 1) around the cores of the two LCSs at roughly $\left(0.5\pm0.25,0.5\right)$,
see \Cref{fig:rdg_anisotropy,fig:Bickley_anisotropy}.
The further away from the structure centers, the more quasi-one-dimensional
the diffusion becomes (volume ratio values close to 0), cf.~also \cite{Thiffeault2001,Thiffeault2003}.
In particular, in \Cref{fig:rdg_diffusivity,fig:Bickley_diffusivity} there are strongly diffusive yellowish
filaments almost enclosing the bluish regions. Second, the direction normal to
the LCS boundaries corresponds to the subdominant diffusion direction, which
is several orders of magnitudes weaker than the dominant
diffusion and therefore significantly slower. To leverage the correspondence between
local geodesic averaging and short-time heat flow, imagine a small geodesic radius
$\varepsilon>0$ and a geodesic ball (in the geometry of mixing) of radius $\varepsilon$
attached to points on some dense grid. With the shape and orientation of these geodesic balls
as described in \cref{sec:anisotropy}, one may guess that the geodesic averaging
operator leaves characteristic functions localized on the LCSs almost invariant, much more
than characteristic functions localized on smaller subsets thereof, or on material sets outside
the LCSs.

\Cref{fig:rdg_density,fig:Bickley_density} demonstrate that the material manifolds
equipped with the respective geometry of mixing may be regarded as
consisting of two and six massive components, respectively, connected by an almost massless background.

In other words, a uniform heat distribution localized close to the isotropic core of the
LCSs diffuses both radially and circularly on comparable
time scales. A uniform heat distribution localized on the whole LCS
will diffuse to the exterior on extremely long time scales and is therefore
expected to be extremely slowly decaying, or, in other words, metastable.
We demonstrate that our expectations built from the diffusion tensor
analysis are indeed satisfied by the heat flow animations provided in the
Supplementary Materials 3 \& 4, and discussed in \cref{sec:heatflow}.

\subsection{Variational Characterization of Eigenvalues}

We continue our study of the geometry of mixing by interpreting
the eigenvalues of the Laplace operator from the variational viewpoint, in light of 
the preceding visualizations of the averaged diffusion tensor field and the density.
This establishes some rather explicit connections between the Lagrangian averaged diffusion
tensor field $\overline{D}$ on the one hand and the topography of low eigenfunctions of the
dynamic Laplacian on the other hand.

According to the Courant--Fischer--Weyl min--max principle the eigenvalues of any
Laplace operator can be characterized as follows: For $k\in\mathbb{N}$ let
$W_{k}=\linhull\{w_{1},\ldots,w_{k}\}\subset L^{2}(M,\mathrm{d}x)$
be the $k$-dimensional subspace spanned by the eigenfunctions corresponding
to the $k$ dominant eigenvalues and $W_{k}^{\bot}$ its orthogonal
complement. Then the $k$th eigenvalue is given by
\[
\lambda_{k} = -\inf_{w\in W_{k-1}^{\bot}}\frac{\int_{M}\bar{g}^{-1}(\mathrm{d}w,\mathrm{d}w)\,\mathrm{d}x}{\int_{M}w^{2}\,\mathrm{d}x}.
\]
The infimum is attained exactly by the eigenfunctions corresponding
to $\lambda_{k}$. For a smooth function $w$, the simplest way to
minimize the Rayleigh quotient is to be non-vanishing and  constant.
On connected manifolds, globally constant functions are captured by
the eigenspace of the zero eigenvalue, and functions in the orthogonal
complement must necessarily (i) have variation and (ii) be sign-indefinite.
From an ``eigenfunction-engineering'' perspective, two questions arise for
dominant eigenfunctions: (i) where to change values
and where to remain (almost) constant; and (ii) if changing values locally
then in which direction the most? Of course, the overall constraint is to
have as little variation as possible.

Clearly, it is generally favorable to have variation in regions where any
variation does not come at a high cost. Pointwise, the maximal cost is
determined by the maximal eigenvalue $\mu_{\max}(\overline{D})$ of the
averaged diffusion tensor, which also dominates the trace of $\overline{D}$.
Thus, one would expect variation in regions with low trace, or, in our terminology
above, with low Lagrangian effective diffusivity. As for the direction of strongest
variation, it is pointwise optimal to have the differential $\mathrm{d}w$ point in the
direction of $v_{\min}(\overline{D})$, i.e., orthogonal to the textured structures shown
in \Cref{fig:rdg_diffusivity,fig:Bickley_diffusivity}.

So far, one may think that it is favorable for a low eigenfunction to take some,
say, positive value in the very LCS center, where effective diffusivity is very low,
and some negative value everywhere else (except for a smooth transition).
Due to the orthogonality constraint $\langle w, 1\rangle$, this would introduce a very steep
gradient in the transition zone, which may outbalance the low cost due to low 
$\mu_{\max}(\overline{D})$. Thus, it is advantageous to push the transition zone
outward, thereby increasing the size of the region with low-variation and positive values,
which allows for a less steep transition between the vortex-like LCS region and its surrounding.
Pushing the transition zone further outward into the high diffusivity region increases the
Rayleigh quotient again, and the low eigenfunctions find an optimal balance in this geometry of mixing.

\subsection{Geometric Heat Flow}\label{sec:heatflow}

\begin{figure}
\centering
\subfloat[]{\label{fig:coherent_initial}\includegraphics[height=.35\textwidth]{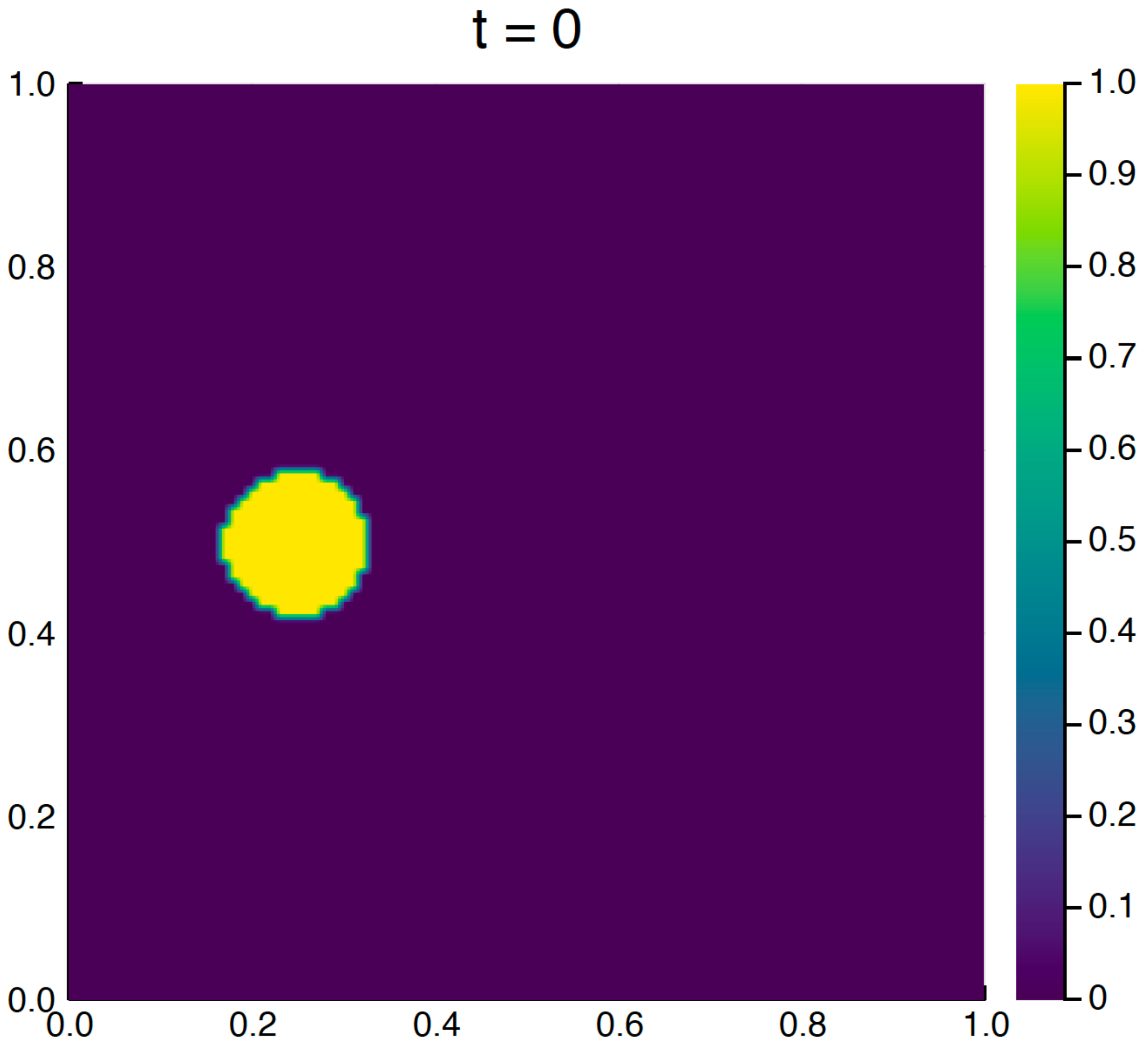}}\qquad
\subfloat[]{\label{fig:coherent_final}\includegraphics[height=.35\textwidth]{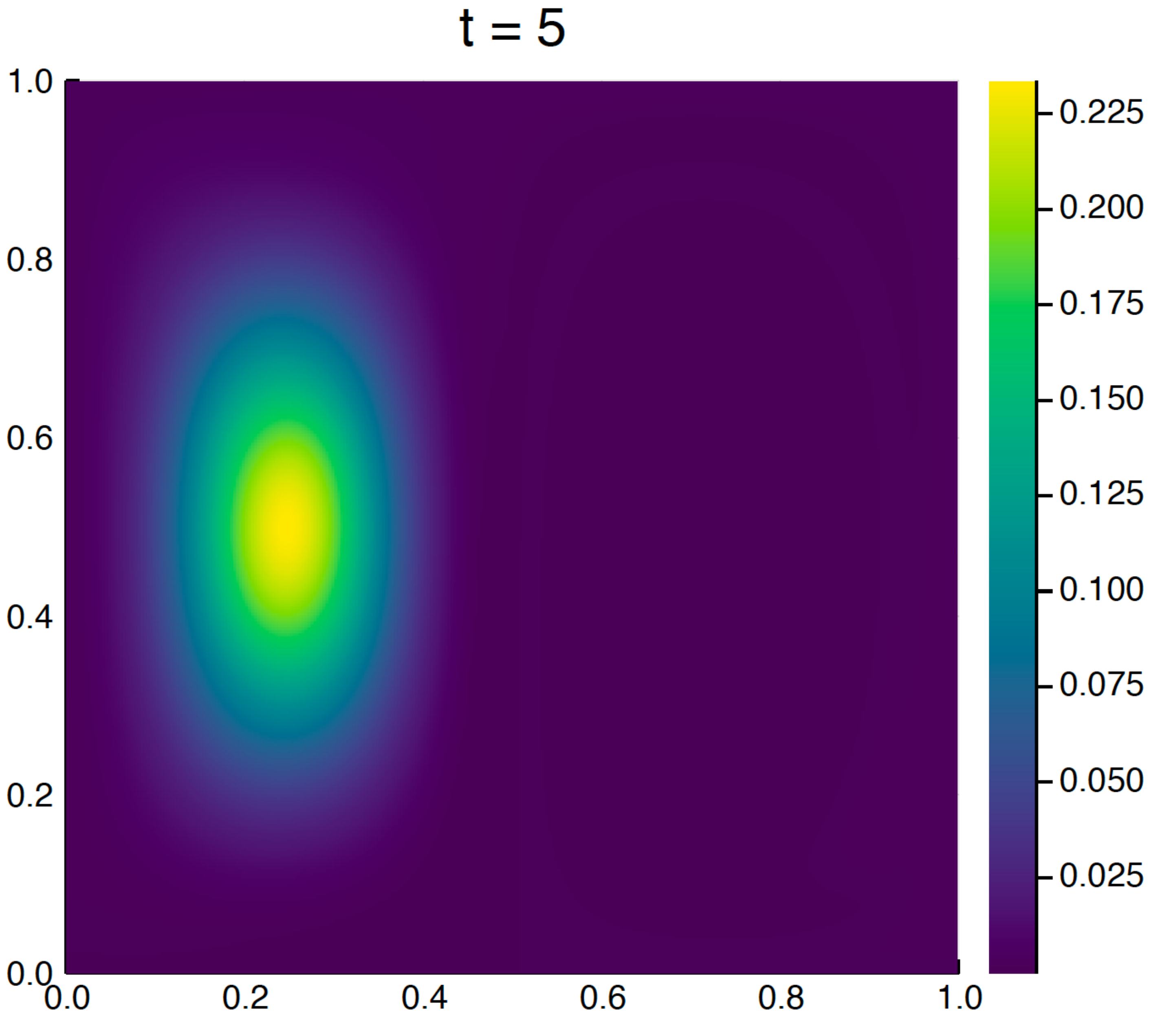}}\\
\subfloat[]{\label{fig:mixing_initial}\includegraphics[height=.35\textwidth]{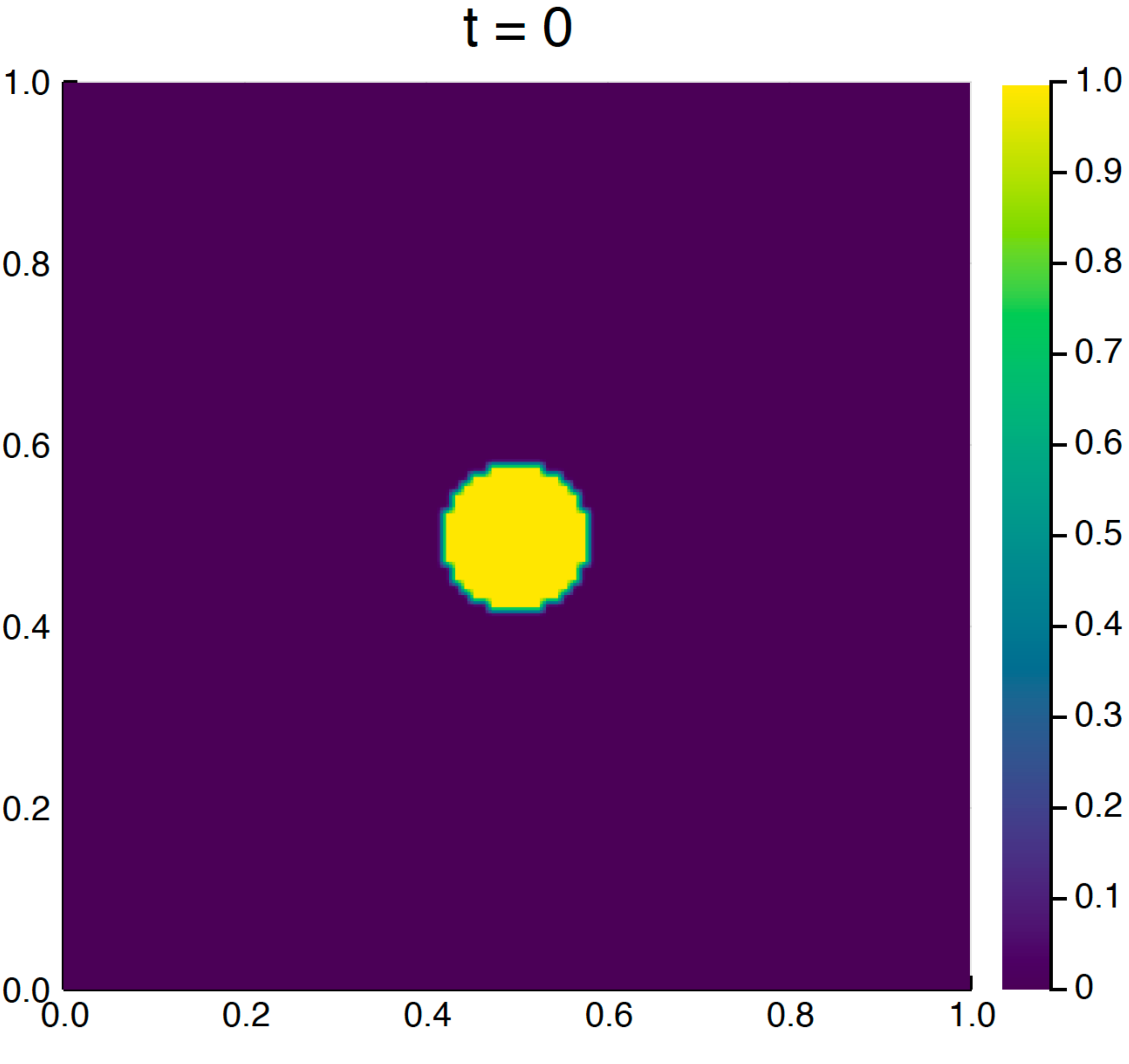}}\qquad
\subfloat[]{\label{fig:mixing_final}\includegraphics[height=.35\textwidth]{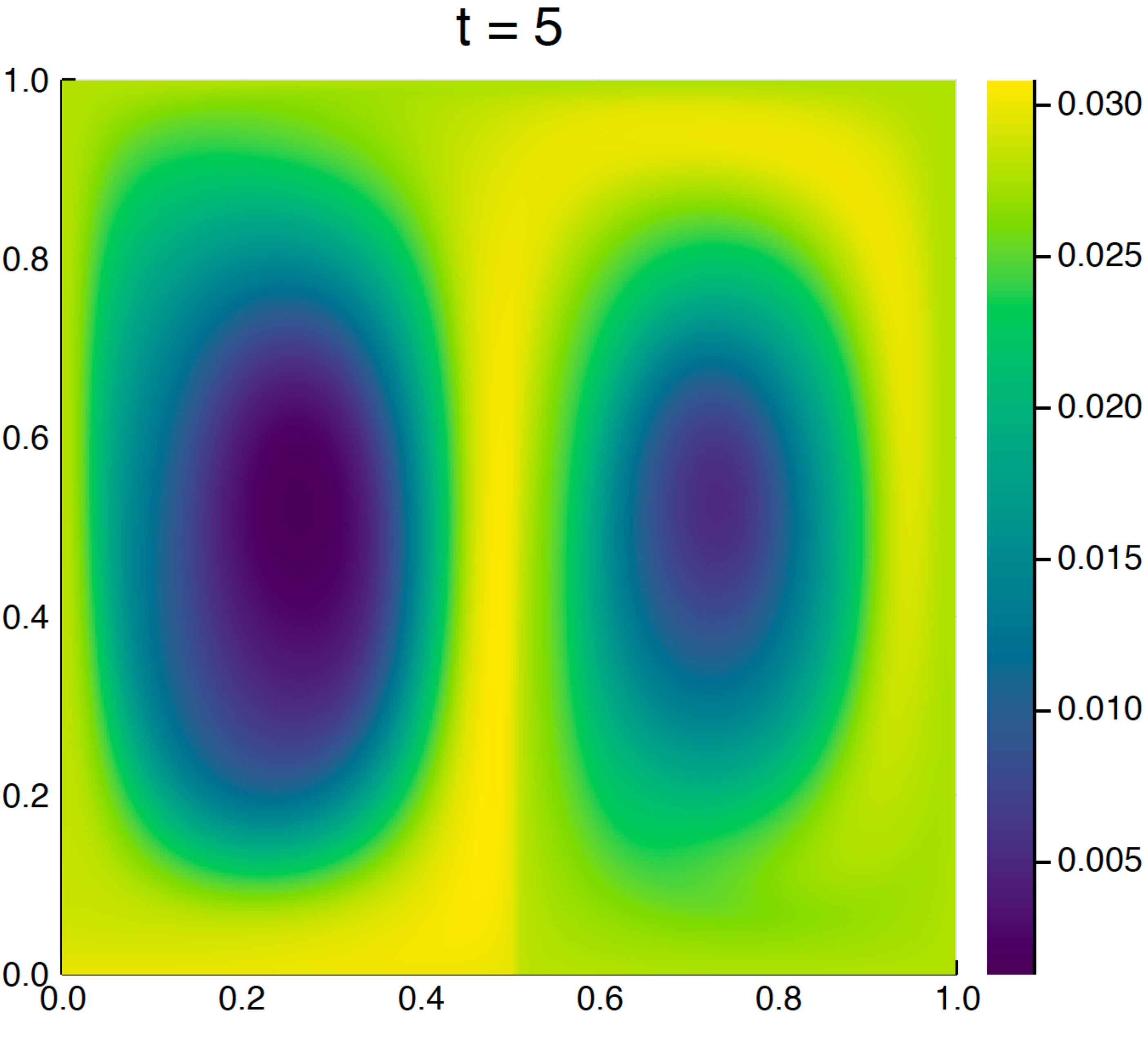}}
\caption{Initial (left) and final (right) scalar densities for an initial condition
localized in the center of an LCS (top) and in the mixing region (bottom).}
\label{fig:heatflow}
\end{figure}

In support of the decomposition of the fluid domain $M$ into regular/coherent
and mixing regions, we look at the action of the geometric heat flow induced by
$\Delta_{\bar{g}}$ on two different initial scalar distributions in the context of
\cref{exa:transient1}, the transient double gyre. To this end, we provide two video
animations of the geometric heat flow in the Supplementary Material 3 and 4.
In the third Supplementary Material, we consider the geometry of mixing and initialize
a scalar quantity in the interior of the left LCS; see \Cref{fig:coherent_initial}.
The scalar quantity is slowly distributed over the LCS, and only very little leaks out;
see \Cref{fig:coherent_final}. 
In the fourth Supplementary Material, we initialize the same amount of some scalar quantity between the
two LCSs in the mixing region; see \Cref{fig:mixing_initial}. Here, the scalar mixes fastly all over the mixing region,
but only very little enters into the two LCSs; see \Cref{fig:mixing_final}. This demonstrates both the mixing 
behavior and the diffusion barrier property of the LCS boundaries.

\subsection{Diffusive Flux Form\label{sec:surface-area-form}}

Another interesting geometric object related to a weighted manifold is its
induced surface area form. Such forms assign a $(d-1)$–volume, which
we will simply refer to as \emph{area}, to $(d-1)$–dimensional parallelepipeds
in tangent space. We restrict our attention to parallelepipeds with
unit $g$–area, whose corresponding $\bar{g}$–area can be interpreted
as the '$\bar{g}$–diffusive flux'. To compute the $\bar{g}$–area
of parallelepipeds of interest we recall a result from linear algebra
\cite[App.~A, Lemma 1]{Froyland2015a}: For an invertible matrix $A\in GL\left(\mathbb{R}^{d}\right)$
and an orthonormal basis $\left(v_{1},\ldots,v_{d}\right)$ one has
\begin{equation*}
\left\lVert A\left(v_{1}\wedge\ldots\wedge v_{d-1}\right)\right\rVert =\det(A)\left\lVert A^{-\top}v_{d}\right\rVert .
\end{equation*}
In our context, given some tangent space $T_{p}M$, $\left(v_{1},\ldots,v_{d}\right)$
shall be an orthonormal basis with respect to $g$. The transformation given by $A\coloneqq\overline{G}^{1/2}$ corresponds to the basis transformation in $T_{p}M$ to Riemannian
normal coordinates, cf.~\cite[Appendix]{Karrasch2015c}. Therefore,
in the new coordinates, $\overline{G}$ takes the canonic Euclidean
form, and area, determinant as well as volume are computed
as in the Euclidean case. In other words, on the left-hand side we
have the $\bar{g}$–area of the parallelepiped spanned by $\left(v_{1},\ldots,v_{d-1}\right)$—our
object of interest— and on the right-hand side we have the density
$\det\left(\overline{G}\right)^{1/2}=\sqrt{\bar{g}}$ discussed
in \cref{sec:density} and the $\bar{g}^{-1}$–norm of the normal
(co-)vector $v_{d}$.

Given a point $p$ in $M$, what is the orientation of a $\left(d-1\right)$-parallelepiped
of unit $g$–area in the tangent space $T_{p}M$ with the least $\bar{g}$–area?
Physically speaking, what is the orientation of a surface element attached to $p$ that admits
the least diffusive flux?
Looking at the left-hand side, we find directly that the parallelepiped
spanned by the eigenvectors of $\overline{G}$ corresponding to the
lowest eigenvalues has minimal $\bar{g}$–area. This is consistent
with the right-hand side in that the normal vector $v_{d}$ is the
eigenvector corresponding to  $v_{\min}\left(\overline{D}\right)$
in that case, and therefore has minimal dual norm.

The area form induced by the geometry of mixing features prominently
in asymptotics of diffusive flux through material surfaces in the vanishing
diffusivity limit \cite{Karrasch2020a}; see also \cite{Haller2018,Haller2020}.

\section{Discussion}
\label{sec:discussion}

\subsection{Generalization to Compressible Flows}
\label{sec:compressible}

The geometric heat-flow approach generalizes by analogy to the setting
when advection is due to a compressible flow. In such contexts, advection--diffusion
of a scalar density $\phi$ as it is carried by a (compressible) fluid with conserved mass
density $\rho$ is modeled by the mass-based equations \cite{Landau1966}, cf.~also \cite{Thiffeault2003},
\begin{subequations}\label{eq:EmADE}
\begin{align}
\frac{D\phi}{Dt}=\partial_{t}\phi+\divergence(\phi V) &=\varepsilon\divergence_{\nu}\, D \mathrm{d}\phi, & \phi(0,\cdot) &= \phi_0, \label{eq:euler_mADE}\\
\frac{D\rho}{Dt}=\partial_t\rho+\divergence(\rho V) &= 0, & \rho(0,\cdot) &= \rho_0. \label{eq:continuity}
\end{align}
\end{subequations}
Here, the (time-dependent) spatial measure $\nu$ corresponds to the fluid's mass
$\mathrm{d}\nu(t)=\rho(t)\mathrm{d}x$, i.e., it has density $\rho$ w.r.t.~the physical volume. \Cref{eq:continuity}
is easily recognized as the continuity equation (or, mass conservation) for $\rho$ and
implies that $\nu(t)$ is the pushforward measure of $\nu_0$ under the flow; cf.~the corresponding
construction in \cite{Froyland2017}.
Moreover, $\divergence_{\nu}\, D \mathrm{d}\phi = \Delta_{g,\nu}$ for the diffusion-adapted
metric $g$ whose coordinate representation is $G=D^{-1}$.

When we represent \cref{eq:EmADE} in Lagrangian coordinates, see \cite{Thiffeault2003}, we obtain
\begin{align*}
\partial_t\varphi &= \varepsilon\Delta_{g(t),\mu(t)}\varphi, & \varphi(0,\cdot) = \phi_0,\\
\partial_t\varrho &= -\varrho\Phi^*(\divergence(V)), & \varrho(0,\cdot) = \rho_0.
\end{align*}
Analogously to the incompressible case, the scalar density $\varphi$ is subject to diffusion
generated by the time-dependent family of Laplace operators $\left(\Delta_{g(t),\mu(t)}\right)_t$,
with $g(t)\coloneqq(\Phi^t)^*g$ the pullback metric, and
$\mathrm{d}\mu(t)=\varrho(t)\mathrm{d}x=(\Phi^t)^*\mathrm{d}\nu(t)$ the pullback measure.
Conservation of mass yields
\begin{equation*}
\mathrm{d}\mu(t)=\mathrm{d}\mu(0)=\rho_0\mathrm{d}x = \mathrm{d}\nu_0,
\end{equation*}
such that---fully analogously to the incompressible case---time dependence enters
the pullback Laplace operators only in the diffusion tensor term. Finally, time averaging
yields the \emph{dynamic Laplacian} operator $\Delta_{\bar{g},\nu_0}$ defined
in \cite{Froyland2017} for compressible flows. In summary, both dynamic Laplacian
constructions given in \cite{Froyland2015a,Froyland2017} can be derived
from a Lagrangian advection--diffusion viewpoint.

\subsection{Spectral Analysis and the Role of Eigenfunctions}
\label{sec:detection}

In our framework, we consider solutions of the self-adjoint elliptic eigenproblem
\begin{equation*}
\Delta_{\bar{g}}w_{n}=\lambda_{n}w_{n},
\end{equation*}
possibly with natural Neumann boundary conditions on $\partial M$. 
As summarized in \cref{sec:metastability}, metastable
sets are identified by separating the dominant part of the spectrum from the rest via
an eigengap. Subsequently, a metastable decomposition is extracted from the associated eigenfunctions.

In our Riemannian manifold framework here, there is an additional spectral interpretation
method that comes into consideration.\footnote{This interpretation is put forward in the recent coherent structure literature: Froyland is alluding to this interpretation in the motivating examples given in \cite{Froyland2015a}, and Banisch \& Koltai \cite{Banisch2017} use eigenfunctions to define \emph{transport coordinates} in the style of \emph{diffusion coordinates} of \cite{Coifman2006}.} It has been developed in the field of
\emph{manifold learning} in the realm of diffusion maps \cite{Belkin2003,Coifman2006}.
The typical problem considered there is the following: given a sample from a \emph{connected}
manifold, that is low-dimensional but embedded in some high-dimensional Euclidean
space, one seeks coordinates that parameterize the manifold intrinsically, without
resorting to the high-dimensional coordinates in the ambient space. This can be
achieved by analyzing the eigenvectors of some graph Laplacian. As discussed in
detail in \cite{Dsilva2015}, one observes the following phenomena as one scans
through the sequence of eigenvectors: the first non-trivial eigenvector [called
\emph{unique eigendirections} in \cite{Dsilva2015}] defines a coordinate direction. It
may be followed either by an eigenvector inducing the same coordinate direction,
however, with a higher frequency of variation (called \emph{repeated eigendirections} 
in \cite{Dsilva2015}), or another unique eigendirection. The purpose of \cite{Dsilva2015}
is to propose an algorithm that distinguishes unique from repeated eigendirections.

The manifold learning problem, however, is different in two aspects to the metastability
problem. First, there is no expectation on smallness of associated eigenvalues nor an
expectation of the existence of a spectral gap in manifold learning. That implies that
one is interested exclusively in consecutive, dominant eigenfunctions in the
metastability analysis \cite{Davies1982}, whereas the coordinate-inducing
eigenfunctions may be scattered \cite{Dsilva2015}. Second, one is interested in
plateaus of eigenfunctions in the metastability analysis, whereas one is interested in
their variation in manifold learning.

These two dichotomous approaches to the interpretation of Laplacian spectra are
potentially reconciled as follows: In the dominant spectrum close to zero, one expects
to find information on the metastable sets, and subsequently coordinate-inducing
eigenfunctions and higher (mixed) harmonics thereof on the almost-disconnected, metastable sets.
\Cref{fig:spectral_decomposition} schematically illustrates the role of different
eigenfunctions in the ideal decoupled and the weakly coupled manifold cases. The
distinction of eigenfunctions and their interpretation in terms of metastable sets or coordinates remains a challenging research question.

\begin{figure}
	\centering
	\subfloat[]{
		\begin{tikzpicture}[scale=1,very thick]
		\draw[-triangle 45] (0.5,0) -- (-6,0);
		
		\draw (0,0) +(0,1mm) -- ++(0,-1mm) node[below]{$0$};
		\node[below=28pt] (O) {components};
		\path (0,0)coordinate(O)  (-2,0)coordinate(G1)   (-5.8,0)coordinate(G2);
		\path plot[mark options={black}, mark=x, mark size=2mm] coordinates { 
			(0,0)
		};
		\path plot[mark options={blue}, mark=x, mark size=1.5mm] coordinates { 
			(-3.5,0) (-4.8,0) (-5.5,0)
		};
		\path plot[mark options={olive}, mark=x, mark size=1.5mm] coordinates { 
			(G1)   (-2.3,0)   (-2.8,0) (-4.0,0)
		};
		
		\draw[red,decorate,decoration={brace,mirror,raise=16pt}] (O) -- (G1) 
		node[pos=0.5,above=20pt]{eigengap};
		\draw[blue] ([yshift=-12pt]G1) +(0,1mm) -- ++(0,-1mm);
		\draw[blue,-triangle 45] ([yshift=-12pt]G1) -- ([yshift=-12pt]G2)
		node[pos=0.5,below]{higher harmonics}
		node[black,pos=0.5, below=15pt]{ \& }
		node[olive,pos=0.5, below=30pt]{coordinates};
		\end{tikzpicture}}
	\subfloat[]{
		\begin{tikzpicture}[scale=1,very thick]
		\draw[-triangle 45] (0.5,0) -- (-6,0);
		
		\draw (0,0) +(0,1mm) -- ++(0,-1mm) node[below]{$0$};
		\node[below=28pt] (O) {component};
		\path (0,0)coordinate(O) (-0.3,0)coordinate(R) (-0.9,0)coordinate(L) (-2,0)coordinate(G1)   (-5.8,0)coordinate(G2);
		\path plot[mark options={black}, mark=x, mark size=1.5mm] coordinates { 
			(0,0)
		};
		\path plot[mark options={blue}, mark=x, mark size=1.5mm] coordinates { 
			(-3.5,0) (-4.8,0) (-5.5,0)
		};
		\path plot[mark options={olive}, mark=x, mark size=1.5mm] coordinates { 
			(G1)   (-2.3,0)   (-2.8,0) (-4.0,0)
		};
		\path plot[mark options={red}, mark=x, mark size=1.5mm] coordinates { 
			(R) (-0.5,0)  (-0.7,0) (L)
		};
		
		\draw[red] ([yshift=12pt]R) +(0,1mm) -- ++(0,-1mm)
		([yshift=12pt]L) +(0,1mm) -- ++(0,-1mm);
		\draw[red] ([yshift=12pt]R) -- ([yshift=12pt]L)
		node[pos=0.5,above=3pt]{quasi-components};
		\draw[blue] ([yshift=-12pt]G1) +(0,1mm) -- ++(0,-1mm);
		\draw[blue,-triangle 45] ([yshift=-12pt]G1) -- ([yshift=-12pt]G2)
		node[pos=0.5,below]{higher harmonics}
		node[black,pos=0.5, below=15pt]{ \& }
		node[olive,pos=0.5, below=30pt]{coordinates};
		\end{tikzpicture}}
	\caption{Spectral structure in the decoupled (a) and perturbed (b) cases. The quasi-components indicate connected coherent sets corresponding to almost invariant sets and 
		arise from the unperturbed zero eigenvalues associated with connected components. 
		The coordinates and higher harmonics correspond to first and higher order Laplace--Beltrami eigenfunctions associated with one or a superposition
		of the components or quasi-components.}
	\label{fig:spectral_decomposition}
\end{figure}

\begin{figure}
	\centering
	\includegraphics{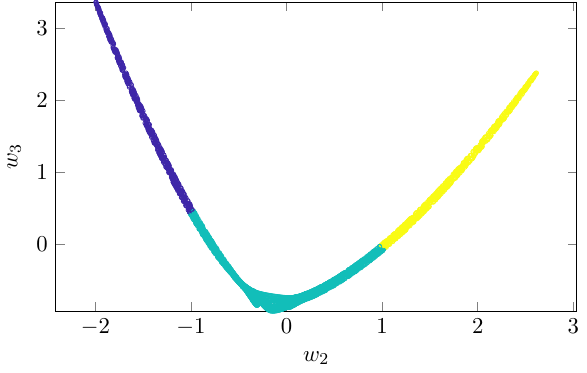}
	\caption{Dominant diffusion coordinates for the rotating double-gyre flow.
		Shown are the second and third eigenfunctions, the flat first eigenfunction
		has been omitted. Each point corresponds to a sample point in the
		flow domain. The coloring corresponds to the coloring of the \textsc{k-means}
		clustering result in \Cref{fig:rdg_clustering}.}
	\label{fig:rdg_coordinates}
\end{figure}

\begin{figure}
\subfloat[$w_{4}$]{\label{fig:rdg_w4}
\begin{tikzpicture}

\begin{axis}[%
width=.33\textwidth,
height=.33\textwidth,
at={(1.269in,0.642in)},
scale only axis,
point meta min=-2.64430776834068,
point meta max=3.793813751093,
axis on top,
xmin=-0.000978473581213307,
xmax=1.00097847358121,
xtick={0,1},
ymin=-0.000978473581213307,
ymax=1.00097847358121,
ytick={0,1},
axis background/.style={fill=white},
colormap={mymap}{[1pt] rgb(0pt)=(0.2081,0.1663,0.5292); rgb(1pt)=(0.211624,0.189781,0.577676); rgb(2pt)=(0.212252,0.213771,0.626971); rgb(3pt)=(0.2081,0.2386,0.677086); rgb(4pt)=(0.195905,0.264457,0.7279); rgb(5pt)=(0.170729,0.291938,0.779248); rgb(6pt)=(0.125271,0.324243,0.830271); rgb(7pt)=(0.0591333,0.359833,0.868333); rgb(8pt)=(0.0116952,0.38751,0.881957); rgb(9pt)=(0.00595714,0.408614,0.882843); rgb(10pt)=(0.0165143,0.4266,0.878633); rgb(11pt)=(0.0328524,0.443043,0.871957); rgb(12pt)=(0.0498143,0.458571,0.864057); rgb(13pt)=(0.0629333,0.47369,0.855438); rgb(14pt)=(0.0722667,0.488667,0.8467); rgb(15pt)=(0.0779429,0.503986,0.838371); rgb(16pt)=(0.0793476,0.520024,0.831181); rgb(17pt)=(0.0749429,0.537543,0.826271); rgb(18pt)=(0.0640571,0.556986,0.823957); rgb(19pt)=(0.0487714,0.577224,0.822829); rgb(20pt)=(0.0343429,0.596581,0.819852); rgb(21pt)=(0.0265,0.6137,0.8135); rgb(22pt)=(0.0238905,0.628662,0.803762); rgb(23pt)=(0.0230905,0.641786,0.791267); rgb(24pt)=(0.0227714,0.653486,0.776757); rgb(25pt)=(0.0266619,0.664195,0.760719); rgb(26pt)=(0.0383714,0.674271,0.743552); rgb(27pt)=(0.0589714,0.683757,0.725386); rgb(28pt)=(0.0843,0.692833,0.706167); rgb(29pt)=(0.113295,0.7015,0.685857); rgb(30pt)=(0.145271,0.709757,0.664629); rgb(31pt)=(0.180133,0.717657,0.642433); rgb(32pt)=(0.217829,0.725043,0.619262); rgb(33pt)=(0.258643,0.731714,0.595429); rgb(34pt)=(0.302171,0.737605,0.571186); rgb(35pt)=(0.348167,0.742433,0.547267); rgb(36pt)=(0.395257,0.7459,0.524443); rgb(37pt)=(0.44201,0.748081,0.503314); rgb(38pt)=(0.487124,0.749062,0.483976); rgb(39pt)=(0.530029,0.749114,0.466114); rgb(40pt)=(0.570857,0.748519,0.44939); rgb(41pt)=(0.609852,0.747314,0.433686); rgb(42pt)=(0.6473,0.7456,0.4188); rgb(43pt)=(0.683419,0.743476,0.404433); rgb(44pt)=(0.71841,0.741133,0.390476); rgb(45pt)=(0.752486,0.7384,0.376814); rgb(46pt)=(0.785843,0.735567,0.363271); rgb(47pt)=(0.818505,0.732733,0.34979); rgb(48pt)=(0.850657,0.7299,0.336029); rgb(49pt)=(0.882433,0.727433,0.3217); rgb(50pt)=(0.913933,0.725786,0.306276); rgb(51pt)=(0.944957,0.726114,0.288643); rgb(52pt)=(0.973895,0.731395,0.266648); rgb(53pt)=(0.993771,0.745457,0.240348); rgb(54pt)=(0.999043,0.765314,0.216414); rgb(55pt)=(0.995533,0.786057,0.196652); rgb(56pt)=(0.988,0.8066,0.179367); rgb(57pt)=(0.978857,0.827143,0.163314); rgb(58pt)=(0.9697,0.848138,0.147452); rgb(59pt)=(0.962586,0.870514,0.1309); rgb(60pt)=(0.958871,0.8949,0.113243); rgb(61pt)=(0.959824,0.921833,0.0948381); rgb(62pt)=(0.9661,0.951443,0.0755333); rgb(63pt)=(0.9763,0.9831,0.0538)},
colorbar
]
\addplot [forget plot] graphics [xmin=-0.000978473581213307,xmax=1.00097847358121,ymin=-0.000978473581213307,ymax=1.00097847358121] {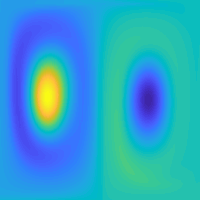};
\end{axis}
\end{tikzpicture}%
}\quad\subfloat[$w_{5}$]{
\begin{tikzpicture}

\begin{axis}[%
width=.33\textwidth,
height=.33\textwidth,
at={(1.269in,0.642in)},
scale only axis,
point meta min=-4.26916488938344,
point meta max=4.47867006432183,
axis on top,
xmin=-0.000978473581213307,
xmax=1.00097847358121,
xtick={0,1},
ymin=-0.000978473581213307,
ymax=1.00097847358121,
ytick={0,1},
axis background/.style={fill=white},
colormap={mymap}{[1pt] rgb(0pt)=(0.2081,0.1663,0.5292); rgb(1pt)=(0.211624,0.189781,0.577676); rgb(2pt)=(0.212252,0.213771,0.626971); rgb(3pt)=(0.2081,0.2386,0.677086); rgb(4pt)=(0.195905,0.264457,0.7279); rgb(5pt)=(0.170729,0.291938,0.779248); rgb(6pt)=(0.125271,0.324243,0.830271); rgb(7pt)=(0.0591333,0.359833,0.868333); rgb(8pt)=(0.0116952,0.38751,0.881957); rgb(9pt)=(0.00595714,0.408614,0.882843); rgb(10pt)=(0.0165143,0.4266,0.878633); rgb(11pt)=(0.0328524,0.443043,0.871957); rgb(12pt)=(0.0498143,0.458571,0.864057); rgb(13pt)=(0.0629333,0.47369,0.855438); rgb(14pt)=(0.0722667,0.488667,0.8467); rgb(15pt)=(0.0779429,0.503986,0.838371); rgb(16pt)=(0.0793476,0.520024,0.831181); rgb(17pt)=(0.0749429,0.537543,0.826271); rgb(18pt)=(0.0640571,0.556986,0.823957); rgb(19pt)=(0.0487714,0.577224,0.822829); rgb(20pt)=(0.0343429,0.596581,0.819852); rgb(21pt)=(0.0265,0.6137,0.8135); rgb(22pt)=(0.0238905,0.628662,0.803762); rgb(23pt)=(0.0230905,0.641786,0.791267); rgb(24pt)=(0.0227714,0.653486,0.776757); rgb(25pt)=(0.0266619,0.664195,0.760719); rgb(26pt)=(0.0383714,0.674271,0.743552); rgb(27pt)=(0.0589714,0.683757,0.725386); rgb(28pt)=(0.0843,0.692833,0.706167); rgb(29pt)=(0.113295,0.7015,0.685857); rgb(30pt)=(0.145271,0.709757,0.664629); rgb(31pt)=(0.180133,0.717657,0.642433); rgb(32pt)=(0.217829,0.725043,0.619262); rgb(33pt)=(0.258643,0.731714,0.595429); rgb(34pt)=(0.302171,0.737605,0.571186); rgb(35pt)=(0.348167,0.742433,0.547267); rgb(36pt)=(0.395257,0.7459,0.524443); rgb(37pt)=(0.44201,0.748081,0.503314); rgb(38pt)=(0.487124,0.749062,0.483976); rgb(39pt)=(0.530029,0.749114,0.466114); rgb(40pt)=(0.570857,0.748519,0.44939); rgb(41pt)=(0.609852,0.747314,0.433686); rgb(42pt)=(0.6473,0.7456,0.4188); rgb(43pt)=(0.683419,0.743476,0.404433); rgb(44pt)=(0.71841,0.741133,0.390476); rgb(45pt)=(0.752486,0.7384,0.376814); rgb(46pt)=(0.785843,0.735567,0.363271); rgb(47pt)=(0.818505,0.732733,0.34979); rgb(48pt)=(0.850657,0.7299,0.336029); rgb(49pt)=(0.882433,0.727433,0.3217); rgb(50pt)=(0.913933,0.725786,0.306276); rgb(51pt)=(0.944957,0.726114,0.288643); rgb(52pt)=(0.973895,0.731395,0.266648); rgb(53pt)=(0.993771,0.745457,0.240348); rgb(54pt)=(0.999043,0.765314,0.216414); rgb(55pt)=(0.995533,0.786057,0.196652); rgb(56pt)=(0.988,0.8066,0.179367); rgb(57pt)=(0.978857,0.827143,0.163314); rgb(58pt)=(0.9697,0.848138,0.147452); rgb(59pt)=(0.962586,0.870514,0.1309); rgb(60pt)=(0.958871,0.8949,0.113243); rgb(61pt)=(0.959824,0.921833,0.0948381); rgb(62pt)=(0.9661,0.951443,0.0755333); rgb(63pt)=(0.9763,0.9831,0.0538)},
colorbar
]
\addplot [forget plot] graphics [xmin=-0.000978473581213307,xmax=1.00097847358121,ymin=-0.000978473581213307,ymax=1.00097847358121] {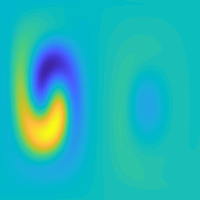};
\end{axis}
\end{tikzpicture}%
}
\caption{Rotating double-gyre flow. Fourth (a) and fifth (b) eigenfunctions, which can be
interpreted as inducing radial (a) and angular (b) coordinates on the left LCS.}
\label{fig:coordinate_eigenfunctions}
\end{figure}

As a demonstration, consider the rotating double-gyre flow, \cref{exa:transient1}.
We used the eigenfunctions $u_2$ and $u_3$ to extract a metastable
decomposition into the two LCSs and the surrounding mixing region in \cref{sec:examples};
see also \Cref{fig:rdg_coordinates} for a typical diffusion--coordinate plot of these.
Eigenfunctions related to eigenvalues following the first two non-trivial ones 
can now be interpreted as inducing the radial and some
``angular'' coordinates on the two LCSs by tracing across the level sets from low to
high values; see \Cref{fig:coordinate_eigenfunctions}.

\subsection{Discretization Aspects}
\label{sec:discretization}

In this section, we collect some thoughts regarding discretization aspects of the
continuous averaged Lagrangian diffusion framework, provide relations to previously
published trajectory-based computational approaches and indicate future research
prospects.

Regarding the discretization of the dynamic Laplacian, there exist two
approaches: (i) start from the continuous formulation and consider different
discretizations or (ii) start from discrete approaches and target for the dynamic
Laplacian in some infinite-sample limit.

As for the first, this refers to the classic and broad field of numerics for partial
differential equations (PDEs). The problem at hand is to discretize the second-order
differential operator $\Delta_{\bar{g}}$, and classic approaches that have been employed
include the finite-difference discretization \cite{Froyland2015a,Froyland2017}, the projection on a 
finite-dimensional subspace of admissible functions spanned by radial basis
functions with compact support \cite{Froyland2015b}, and the finite-element method (FEM)
\cite{Froyland2018}. In most practical cases, one will work on a domain with 
boundary and hence has to ensure that (homogeneous) boundary conditions are 
satisfied. Moreover, one would want to preserve a crucial feature of
$\Delta_{\bar{g}}$, namely its self-adjointness. Out of the afore-mentioned discretization 
schemes, only the FEM approach \cite{Froyland2018} manages to fulfill both these 
requirements exactly and in addition admits an implementation suitable to handle
sparse and incomplete data.

As for the second, this leads to the well-established field of diffusion maps.
The central result is, roughly speaking, that a graph Laplacian matrix built from
a sample from some Riemannian manifold converges toward the classic
Laplace operator on the manifold. There is an abundance of varieties here:
which graph Laplacian (normalized or unnormalized) to consider, which kernel
function to use (fixed- or variable-bandwidth) and which manifolds (with or 
without boundary, compact or not) to allow, and finally, what kind of convergence
to obtain \cite{Coifman2006,Luxburg2008,Berry2016}.

Briefly, let $(x_i)_i$ denote the samples from the manifold $M$, then define the
entries $w_{ij}$ of the weight matrix $W$ associated with the pair $(x_i,x_j)$ via a
kernel $k\colon\mathbb{R}_{\geq 0}\to\mathbb{R}_{\geq 0}$ and the distance
$\dist(x_i,x_j)$, i.e., $w_{ij}=k(\dist(x_i,x_j))$, where $k$ may additionally depend on
$x_i$ and $x_j$. A classic choice is given by the Gaussian
$k(w) = \exp(-w^2/\sigma^2)$, but the exact design of the kernel turns out to be a
powerful tool on its own \cite{Berry2016a}. From the weight matrix $W = (w_{ij})$, one
constructs the diagonal \emph{degree matrix} $D$, with diagonal entries
$d_i = \sum_j w_{ij}$, and finally $I-D^{-1}W$ is the normalized graph Laplacian, an
approximation to the negative Laplace operator. This goes back to the normalized
spectral clustering in \cite{Shi2000} and admits a plethora of modifications and
extensions \cite{Coifman2006}.

In this spirit, several suggestions from the recent literature admit an interpretation as
approximations of averaged Lagrangian diffusion. Since the dynamic Laplacian is
defined as the average of pullback Laplace operators, it is natural to consider its
approximation via the average of graph Laplacian approximations of the pullback
Laplace operators \cite{Banisch2017}. This has the advantage that one can use
distances between trajectories at different times (usually Euclidean distances) for the
construction of the graph Laplacians, but has the drawback that the final matrix
representation of $\Delta_{\bar{g}}$ \cite[$\mathbf{Q}_{\varepsilon}$ in Eq.~(19)]{Banisch2017} does not preserve self-adjointness \cite[p.~8, comment (a)]{Banisch2017}.
However, it approximates the pointwise action of $\Delta_{\bar{g}}$ on functions
in the infinite-sample limit \cite[Thm.~3]{Banisch2017}.

A---challenging---alternative might be to operate directly on the intrinsic geometry induced by
$\bar{g}$. In the spirit of graph Laplacians, this requires the
computation/approximation of geodesic distances in the $\bar{g}$-metric. If the
tensor field can be evaluated on a rather fine grid, this would allow the computation
of geodesic distances by solving the (anisotropic) eikonal equation. For this purpose,
very efficient algorithms such as the fast-marching algorithm or, more generally level
set methods, are, in principle, available. Alternatively, one may approximate geodesic
distances to the immediate grid neighbors from the metric tensors computed on the
grid and extend the distance metric to the whole grid by computing lengths of
shortest paths.

Another line of research has suggested to define notions of distance, or, even
weaker, of adjacency (in a graph sense) between trajectories; cf.~\cite{Froyland2015}
for an approach that introduces a metric on the extended state space, which is
then applied to compute distances between trajectories and cluster centers as
required by the k-means algorithm. On the one hand,
\cite{Hadjighasem2016} proposed to use time-averages of pairwise distances,
\emph{dynamical distances}, as $\dist$-function in the graph Laplacian construction. 
Thereby, pairs of trajectories with large dynamical distance get a low weight which makes it
unlikely that they are diagnosed as members of the same cluster/LCS. On the other
hand, \cite{Padberg-Gehle2017} proposed a complementary approach: regard two
trajectories as adjacent in a graph whose nodes correspond to trajectories if they ever get closer than a
specified distance threshold. Roughly speaking, while \cite{Hadjighasem2016}
specifically punishes large dynamic distances, \cite{Padberg-Gehle2017} rewards close
distances. Notably, the degree field derived from the graph adjacency matrix in
\cite{Padberg-Gehle2017} has been independently proposed as a measure for mixing
potential called trajectory encounter volume \cite{Rypina2017}. As discussed earlier,
one may interpret this equally as a measure of Lagrangian effective diffusivity.
While this correspondence to physical quantities is appealing,
the graph constructed in \cite{Padberg-Gehle2017} is merely an unweighted and undirected
graph. For mathematical analysis, more structure such as an underlying metric or 
even a Riemannian geometry generating adjacency would be desirable. In any case,
for the methods discussed in this paragraph it remains unclear what a continuous
infinite-sample limit could look like and whether it carries a nice mathematical 
structure. Ignoring such convergence aspects, there remains a high degree of analogy
between the dynamic Laplacian and its associated geodesic distance in the geometry
of mixing on the one hand, and ad hoc dynamic distances which enter classic Laplace
operator discretizations on the other hand; see also \cite[Sect.~3.4]{Froyland2018} for
an interpretation of an FEM discretization of the dynamic Laplacian as a graph Laplacian matrix.

One common feature of the above-listed methods is that they embed each trajectory
in some neighborhood via distances, taking into account the dynamics.
This is in contrast to an independent set of coherent structure detection
methods, for instance \cite{Mezic2010,Mancho2013,Mundel2014,Haller2016,Fabregat2016}. These methods
retrieve information from time averages of observables along trajectories, which are
viewed individually, without an immediate neighborhood relationship. The expectation 
then is that coherent structures reveal themselves as sets of material points showing
similar statistics within the structure, and different statistics compared to the exterior.
In other words, a neighborhood relationship is built only after observing trajectory 
statistics.

\subsection{Connections to Geodesic LCS Approaches}
\label{sec:Haller_methods}

There is another group of methods for finding boundaries of coherent
structures in purely advective flows, developed by Haller and co-workers
\cite{Haller2015,Haller2013a,Karrasch2015,Farazmand2014a}. These
build on global variational principles formulated in terms of the stretching or
shear of material boundaries. These approaches naturally involve the
Cauchy–Green (CG) strain tensor field $C\coloneqq D\Phi(T)^{\top}D\Phi(T)$,
which we interpret here as the pullback metric $\Phi(T)^{*}\mathfrak{g}$, see
\Cref{rem:pullbackmetric}. These methods evaluate the dynamics
at two time instances, an initial $t=0$ and a final one $t=T$. Earlier
LCS methods have only used the logarithm of the maximal
eigenvalue of $C$, well known as the \emph{finite-time Lyapunov exponent
(FTLE)}, for visual inference of coherent structures.

First, observe the following tight relation between the Cauchy–Green
strain tensor $C$ and the two-point averaged, diffusion-adapted metric
tensor $\overline{G}$ in physical $\mathfrak{g}$–normal coordinates at
some point $p\in M$. Then $\overline{G}$ has the coordinate representation
$\overline{G}=2\left(G(0)^{-1}+G(T)^{-1}\right)^{-1}=2\left(I+C^{-1}\right)^{-1}$,
where $I$ is the identity matrix. Clearly, $C$ has eigenvalues $\mu_{\min}=\mu_{\min}(C)\leq\mu_{\max}(C)=\mu_{\max}$
with eigenvectors\footnote{Haller and co-workers usually employ the notation $(\lambda_{i},\xi_{i})$
for eigenpairs.} $v_{\min}(C)$ and $v_{\max}(C)$ if and only if $\overline{D}=\overline{G}^{-1}$
has eigenvalues $\mu_{\min}\left(\overline{D}\right)=\frac{1}{2}\left(1+\frac{1}{\mu_{\max}}\right)\leq\frac{1}{2}\left(1+\frac{1}{\mu_{\min}}\right)=\mu_{\max}\left(\overline{D}\right)$
with eigenvectors $v_{\min}\left(\overline{D}\right)=v_{\max}(C)$
and $v_{\max}\left(\overline{D}\right)=v_{\min}(C)$. In other
words, the minor CG–eigendirection $v_{\min}(C)$ corresponds to the
dominant $\overline{G}$–diffusion direction.
Moreover, in the volume-preserving case, one has $\mu_{\min}=\mu_{\max}^{-1}$,
and therefore the anisotropy ratio for $\overline{D}$ 
\[
\frac{1+\mu_{\max}}{1+\frac{1}{\mu_{\max}}}=\mu_{\max}
\]
is equal to the dominant CG-eigenvalue. Thus, the logarithm
of the anisotropy ratio shown in the figures in \cref{sec:tensor-field}
corresponds to an accordingly defined multiple-time-step FTLE up to rescaling.

Next, let us briefly recall the variational formulations for elliptic
(coherent vortex boundaries) and parabolic (jet cores) Lagrangian
coherent structures (LCSs) in two-dimensional flows, using our notation.
Following \cite{Haller2013a}, boundaries of elliptic LCSs are sought
as the outermost \emph{closed} stationary curves of the averaged strain
functional
\[
Q(\gamma)=\int_{0}^{b}\frac{\lvert r'(s)\rvert_{\Phi(T)^{*}g}}{\lvert r'(s)\rvert_{g}}\mathrm{d}s,
\]
where $r$ is a parameterization of a material curve $\gamma\subset M$.
The integrand compares pointwise the magnitude of  curve velocity
$r'(s)$ after push forward by $D\Phi$ with its original magnitude.
Equivalently, by equipping $M$ with the pullback metric $C=\Phi(T)^{*}g$,
the domain $M$ is geometrically deformed, in principle as we do in
our approach here with $\bar{g}$. In these terms, the length of $r'(s)$
in the deformed geometry $(M,C)$ is compared with its length in the
original geometry $(M,g)$. By applying Noether's theorem, one obtains
that stationary curves necessarily obey a conservation law, which
corresponds exactly to the integrand, i.e.,
\[
\frac{\lvert r'(s)\rvert_{\Phi(T)^{*}g}}{\lvert r'(s)\rvert_{g}}=\frac{\lvert r'(s)\rvert_{g(T)}}{\lvert r'(s)\rvert_{g(0)}}=\lambda=\text{const}.,
\]
from which one may deduce tangent line fields $\eta_{\lambda}$. Among
the orbits of these line fields one looks for closed ones. Closed
orbits typically come in continuous one-parameter families, out of
which one picks the outermost. Analogously to index theory for vector
fields, one may employ index theory for line fields to deduce that
any closed orbit of a piecewise differentiable line field must necessarily
enclose at least two singularities\footnote{Singularities of line fields are points
at which the line field is not continuously defined, see
\cite[Chap.~4, Addendum 2]{Spivak1999} and \cite{Delmarcelle1994}. Wedge-type
singularities are characterized by a sector of integral curves running into the
singularity, complemented by a sector of integral curves flowing around the singularity.} of the line fields,
and all enclosed
singularities obey a topological rule, see \cite{Karrasch2015} for
more details. In most cases, relevant closed orbits have been found
to enclose exactly two singularities of wedge type \cite{Karrasch2015}.
These singularities
can be visualized (and numerically detected) by phase portraits of
the eigendirections of the pullback metric $C$, or, analogously, by the dominant
diffusion direction field of $\overline{D}$ as in \cref{sec:tensor-field}.

In practice, the outermost closed stationary curve travels through
regions of high FTLE/anisotropy values. There, the tangent direction
field is almost collinear with $v_{\min}(C)$. From the discussion
in \cref{sec:anisotropy,sec:examples} we conclude that such closed curves
are pointwise very close to the optimal direction for blocking $\bar{g}$–diffusion.
Their deviation from the optimal direction is not very costly in terms
of diffusive flux, but still allows them to close up smoothly under
the conditions of the variational principle. Our ad-hoc considerations here
have been formalized and made rigorous in recent work by Haller \emph{et al.}
\cite{Haller2018,Haller2020}, who provide a variational theory for detecting
material transport barriers to diffusive transport in the original, time-dependent
advection--diffusion process, \cref{eq:euler_ADE,eq:lagr_ADE}.

\begin{figure}
\centering
\subfloat[\label{fig:rdg_lambdakmeans}]{\includegraphics{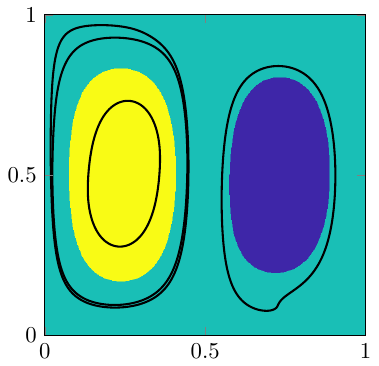}}%
\subfloat[\label{fig:rdg_lambda}]{%
\includegraphics{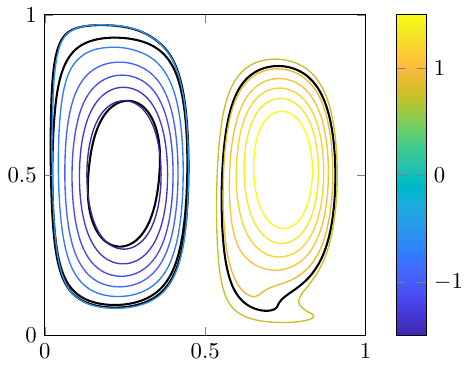}}
\caption{Rotating double-gyre flow (evaluated at initial
and final time points only). Closed $\eta_{\lambda}$–lines (black) on top of
(a) the final clustering result from the $\Delta_{\bar{g}}$–analysis,
and (b) the second eigenfunction $w_{2}$ of $\Delta_{\bar{g}}$. }
\label{fig:Laplace_CG}
\end{figure}

For a simple numerical test, we have overlaid the final clustering
result as well as the second eigenfunction of $\Delta_{\bar{g}}$  for a two
time point-approximation of $\bar{g}$ with all closed $\eta_{\lambda}$–orbits in
\Cref{fig:Laplace_CG}. We find that the eigenvector clustering procedure yields
smaller structures; see \Cref{fig:rdg_lambdakmeans}. A close inspection of the 
level sets of the second eigenfunction, however, reveals a nice, albeit not perfect, visual matching of some level sets with the geodesic LCS results; see \Cref{fig:rdg_lambda}.

In summary, its dimension-independent formulation and its high degree
of consistency with the two-dimensional variational principles suggest
our methodology as a natural extension of these approaches to higher
dimensions. It has proven to be notoriously challenging to extend the variational
ideas to three dimensions by restricting oneself to variational principles
on curves \cite{Blazevski2014,Oettinger2016}.

\subsection{Applications to Geophysical Fluid Dynamics}\label{sec:gfd}

In this section, we compare our methodology to Nakamura's \emph{effective diffusivity}
framework \cite{Nakamura1996,Shuckburgh2003}, which is widely used in the
geophysical fluid dynamics community.

Both methods consider advection--diffusion processes in possibly turbulent fluid flows,
i.e., the advection--diffusion equation (ADE) in the physical domain with a constant diffusion coefficient $\kappa$.\footnote{The assumption of constant diffusion coefficient is delicate in a differential geometry context.
We view it as an assumption on the chosen coordinate system.} In a second step,
the ADE is transformed to a different set of coordinates: In \cite{Nakamura1996}, an
almost Lagrangian coordinate system based on the area inside concentration level
sets of a given tracer density $\phi$ of interest is constructed.
This transformation is considered as advantageous, for it separates ``the reversible
effects of advection from the irreversible effects of advection and diffusion acting
together'' \cite{Shuckburgh2003}. In our context, this is achieved by passing
to Lagrangian coordinates as in \cite{Press1981,Thiffeault2003,Fyrillas2007}---and
therefore literally to a tracer-based coordinate system---which factors out the
advective motion, but keeps the joint action of advection and diffusion through
tracking deformation and its effect on diffusion in the form of the pullback metric.

In the Nakamura framework, the coordinates are then given by the one-dimensional area coordinate
$A$ and $(d-1)$-dimensional coordinates on concentration level sets.
Clearly, the instantaneous local action of diffusion on the concentration is in the $A$-direction only, 
since there is no $\phi$-gradient along the level set coordinates. Averaging over contours
allows to reduce the $d$-dimensional advection--diffusion equation to
a one-dimensional (along the area coordinate) pure diffusion equation \cite{Nakamura1996}:
\begin{equation}
\partial_t\phi = \partial_A(K_{\textnormal{eff}}(A,t)\partial_A\phi),\label{eq:effective_diffusion}
\end{equation}
and the scalar $K_\textnormal{eff}$ is coined \emph{effective diffusivity}.

In our context, the Eulerian ADE is turned into a pure diffusion equation as well, \cref{eq:lagr_ADE},
however, of full dimension $d$ again. Moreover, we arrive at a diffusion \emph{tensor} induced by the pullback metric.
We comment on the reason and the benefit of this increased complexity later.

The next important concept in Nakamura's framework is that of \emph{equivalent length}
$L_{\textnormal{eq}}$ (in 2D) and \emph{equivalent area} $A_{\textnormal{eq}}$ (in 3D).
For simplicity, we focus on $L_{\textnormal{eq}}$ henceforth, but all arguments apply
analogously to higher-dimensional flows. And so, the equivalent length of a concentration
level set corresponding to the area coordinate $A$ can be obtained from a comparison of
the spatial scalar diffusivity $\kappa$ (from the Eulerian ADE) and the effective diffusivity
$K_{\textnormal{eff}}$ (from \cref{eq:effective_diffusion}), cf.\ \cite[eq.~(6)]{Shuckburgh2003}:
\begin{equation}
K_{\textnormal{eff}}(A,t) = \kappa L_{\textnormal{eq}}^2(A,t), \qquad\text{or equivalently}\qquad \frac{K_{\textnormal{eff}}(A,t)}{L_{\textnormal{eq}}^2(A,t)} = \kappa.\label{eq:equilength}
\end{equation}
It can be shown that $L_{\textnormal{eq}}^2$ is at least (the square of) the physical length of the corresponding concentration level set.

\Cref{eq:equilength} can be interpreted as answering the following question: How would
one have to rescale units or, more generally, deform the local geometry in order to
observe---in the new units/deformed geometry---the original diffusivity $\kappa$?
Conceptually, this procedure is completely analogous to the duality of the diffusion
tensor and the induced metric tensor (with respect to which anisotropic diffusion
is isotropic; cf.~\cref{sec:Euler_ADE}), but also to the perspective that we took 
throughout \cref{sec:tensor-field}: How do length, diffusivity, and volume in the geometry
of mixing relate to the corresponding entities in the original spatial geometry?

How do we see in the geometry of mixing whether diffusion is effectively enhanced or suppressed
(relative to the pure spatial diffusion) via its interaction with advection?
As an early indicator of mixing efficiency, Nakamura \cite{Nakamura1996} 
suggested to look at the equivalent length, for a large $L_\textnormal{eq}$
implies a large effective diffusivity, which is then related to the \emph{mixing region}.
By analogy, the Lagrangian effective diffusivity $1/d\cdot\trace\left(\overline{D}\right)$,
$d$ the dimension of space, plays the role of $K_{\textnormal{eff}}$ and
serves as an indicator for the mixing region.

The effective diffusivity framework is tied to a specific tracer concentration field. For instance,
the fact that diffusion is one-dimensional and the orientation of that single coordinate in physical
space is determined by the concrete concentration field at hand. A concentration field with a different initial
level set topology may yield different results. To obtain physically relevant mixing information,
it is assumed that a ``suitable tracer field'' is considered; cf.\ the effort taken in \cite{Nakamura1996,Shuckburgh2003}
to generate those. Suitability there means, roughly speaking, that its level set topology is already 
reasonably ``equilibrated'' w.r.t.\ the mixing geometry induced by the flow, but at the
same time has sufficiently strong gradients to allow for a computationally robust
transformation to area coordinates. In other words, concentration gradients are roughly
aligned with the direction of slowest diffusion, which is of greatest interest.

In contrast, our geometry of mixing provides information about the
``mixing ability'' \cite{Shuckburgh2003} or ``mixing potential'' \cite{Rypina2017} of
the flow, independently from a concrete concentration field.
The Lagrangian ADE remains of full dimension and invokes an
effective diffusion tensor to account for all possible level set topologies. Moreover,
this diffusion tensor field and its corresponding Laplace operator are obtained
computationally from purely advective ODE simulations, there is no need to solve the
advection--diffusion PDE.

\subsection*{Acknowledgments}

The Authors are grateful to Ulrich Pinkall for a useful discussion which triggered this
collaboration. D.K.~would like to thank the organizers of the workshop on
\href{http://www.birs.ca/events/2017/5-day-workshops/17w5048}{``Transport in Unsteady
Flows: from Deterministic Structures to Stochastic Models and Back Again''}
at the BIRS in Banff, Canada, for the invitation and an opportunity to present
parts of this work. Very inspiring discussions there have culminated in
\cref{sec:gfd}. D.K.~would also like to thank George Haller for helpful discussions.
We wish to thank Christian Ludwig and Georg Wechslberger for technical help with
the TikZ figures, and three anonymous referees, whose thoughtful comments have
helped to significantly improve the presentation in the paper.

For the direction field visualization we used a line integral convolution (LIC) code 
which is based on a Matlab-interface by Sonja Rank and Folkmar Bornemann for
\texttt{C}-code developed by Tobias Preusser \cite{Preusser1999}.
For the eigendecomposition of $\Delta_{\bar{g}}$ and the heat flow simulation we
used Julia code developed by Nate Schilling, based on \cite{Froyland2018}, as part of \href{https://github.com/CoherentStructures/CoherentStructures.jl}{\texttt{CoherentStructures.jl}}. 

D.K.~acknowledges partial support by the Priority Programme SPP 1881 ``Turbulent Superstructures''
of the Deutsche Forschungsgemeinschaft (DFG).
D.K.~and J.K.~acknowledge partial support by the German Research Foundation (DFG), 
Collaborative Research Center SFB-TRR 109 ``Discretization in Geometry and Dynamics''.


\end{document}